\documentclass{amsart}
\usepackage{amssymb,amscd,verbatim}
\numberwithin{equation}{section}

\newtheorem{prop}{Proposition}
\newtheorem{thm}[prop]{Theorem}
\newtheorem{cor}[prop]{Corollary}
\newtheorem{lem}[prop]{Lemma}

\theoremstyle{definition}
\newtheorem{ex}[prop]{Example}
\newtheorem{rem}[prop]{Remark}

\begin{document}

\title[Littlewood-Richardson Cyclage Poset]{A Cyclage Poset Structure \\
for Littlewood-Richardson tableaux}

\author{Mark Shimozono}
\address{Dept.\ of Mathematics \\ Virginia Tech \\
Blacksburg, VA}
\email{mshimo@math.vt.edu}

\begin{abstract} A graded poset structure is defined for the
sets of Littlewood-Richardson (LR) tableaux that count the
multiplicity of an irreducible $gl(n)$-module in the tensor
product of irreducible $gl(n)$-modules corresponding to
rectangular partitions.  This poset generalizes the
cyclage poset on column-strict tableaux defined by Lascoux and
Sch\"utzenberger, and its grading function generalizes the charge
statistic.  It is shown that the polynomials obtained by enumerating
LR tableaux by shape and the generalized charge, are none other than
the Poincar\'e polynomials of isotypic components of the
certain modules supported in the closure of a nilpotent conjugacy class.
In particular explicit tableau formulas are obtained for
the special cases of these Poincar\'e polynomials given by
Kostka-Foulkes polynomials, the coefficient polynomials of two-column
Macdonald-Kostka polynomials, and the Poincar\'e polynomials of isotypic
components of coordinate rings of closures of conjugacy classes of
nilpotent matrices.  These $q$-analogues conjecturally coincide with
$q$-analogues of the number of certain sets of rigged configurations and the
$q$-analogues of LR coefficients defined by the spin-weight generating
functions of ribbon tableaux of Lascoux, Leclerc, and Thibon.
\end{abstract}

\maketitle

\newcommand{\charge}{\mathrm{charge}}
\newcommand{\chargediff}{\Delta\mathrm{c}}

\newcommand{\la}{\lambda}
\newcommand{\wh}{\widehat{w}}\
\newcommand{\Rhat}{\widehat{R}}
\newcommand{\K}{K}
\newcommand{\tK}{\widetilde{K}}
\newcommand{\Ph}{\widehat{P}}
\newcommand{\Qh}{\widehat{Q}}
\newcommand{\Th}{\widehat{T}}
\newcommand{\Vh}{\widehat{V}}
\newcommand{\NE}{\mathrm{NE}}
\newcommand{\SW}{\mathrm{SW}}
\newcommand{\RRR}{\mathcal{R}}
\newcommand{\NN}{\mathbb{N}}
\newcommand{\ZZ}{\mathbb{Z}}
\newcommand{\C}{\mathbb{C}}
\newcommand{\LRC}{c}
\newcommand{\LRT}{\mathrm{LRT}}
\newcommand{\ST}{\mathrm{ST}}
\newcommand{\word}{\mathrm{word}}
\newcommand{\cword}{\mathrm{cword}}
\newcommand{\key}{\mathrm{key}}
\newcommand{\NC}{\mathcal{N}}
\newcommand{\SSS}{\mathcal{S}}
\newcommand{\TTT}{\mathcal{T}}
\newcommand{\LLT}{c}
\newcommand{\sign}{\mathrm{sign}}
\newcommand{\weight}{\mathrm{weight}}
\newcommand{\Roots}{\mathrm{Roots}}
\newcommand{\inner}[2]{\langle #1\,,\,#2\rangle}
\newcommand{\sq}{\times}
\newcommand{\Knuth}{\sim_K}
\newcommand{\pr}{\mathrm{pr}}
\newcommand{\rev}{\mathrm{rev}}

\section{Introduction}

In a series of papers \cite{La} \cite{LS1} \cite{LS2}
Lascoux and Sch\"utzenberger developed the deep theory
of the cyclage poset on column-strict tableaux, to give
combinatorial explanations of properties of the Kostka-Foulkes polynomials,
which are $q$-analogues of weight multiplicities in type A.
In particular they assign to each tableau a nonnegative
integer called the charge and show that the Kostka-Foulkes polynomials
are the $q$-enumeration of column-strict tableaux by the charge
\cite{B} \cite{LS1}.

These polynomials occur as instances of Poincar\'e polynomials
$\K_{\la;R}(q)$ of isotypic components of coordinate rings of
closures of nilpotent conjugacy classes of matrices twisted by
line bundles.  The polynomials
$\K_{\la;R}(q)$ are $q$-analogues of the tensor product multiplicities
given by Littlewood-Richardson coefficients of the form
\begin{equation} \label{LR mult}
  c^R_\la = \inner{s_\la}{s_{R_1} s_{R_2} \dotsb s_{R_t}}
\end{equation}
where $s_\la$ is the Schur polynomial and $R=(R_1,R_2,\dotsc,R_t)$
is a sequence of rectangular partitions.  This multiplicity has
a well-known description as the cardinality of a set of
Littlewood-Richardson tableaux $\LRT(\la;R)$.  Many properties of
Kostka-Foulkes polynomials
have generalizations for the Poincar\'e polynomials $\K_{\la;R}(q)$
\cite{SW}.  In \cite{KS} \cite{SW} many combinatorial conjectures were
proposed to explain these properties.

The key construction of this paper is a direct definition of a
poset structure on Littlewood-Richardson tableaux that generalizes the
cyclage poset on column-strict tableaux.
This new poset is graded by a function $\charge_R$ that generalizes
the charge.  The most important consequence of this construction is a
proof that the Poincar\'e polynomial $\K_{\la;R}(q)$ is given by
the $q$-enumeration of the set $\LRT(\la;R)$ of LR tableaux by $\charge_R$.
Other consequences include proofs of monotonicity
and symmetry properties of the polynomials $\K_{\la;R}(q)$
\cite{KS}; these proofs will appear elsewhere.

Conjecturally the polynomials $\K_{\la;R}(q)$
coincide with the $q$-enumeration of rigged configurations \cite{KS}.
A. N. Kirillov has given a bijection from LR tableaux to rigged
configurations, but the conjecture that it
preserves the appropriate statistics, remains open.

The polynomials $\K_{\la;R}(q)$ also conjecturally coincide with a subfamily
of the $q$-analogues of LR coefficients arising from the spin-weight
generating function over ribbon tableaux \cite{LLT}.  At this time,
little is known about this larger family of polynomials.

While preparing this paper the author discovered that
the polynomials $\K_{\la;R}(q)$ appear as multiplicities
of Schur functions in certain Demazure characters of affine type A,
generalizing a level-one formula of \cite{KMOTU} to arbitrary level.
This and the many connections between the affine type A
crystal theory and the tableau combinatorics in
this paper, will appear elsewhere.  

The first section reviews the definitions of the Poincar\'e polynomials,
the LR tableaux, the action of the symmetric group on LR tableaux
by the generalized automorphisms of conjugation,
the generalized charge statistic, and the main result.
The main construction, the $R$-cocyclage poset structure on LR tableaux,
is introduced in Section \ref{cyclage sec}, where the important features
of this graded poset and its grading function $\charge_R$
are described.  Section \ref{proof outline} sketches the
proof of the main theorem.  The proofs appear in the last two sections.

Thanks to M. Okado for pointing out the reference
\cite{ScWa} which has considerable overlap with this paper and \cite{Sh2}.

\section{Definitions}
This section reviews the definitions necessary to state the main result.

\subsection{The Poincar\'e polynomials $\K_{\la;R}$}
We recall the definition of the Poincar\'e polynomials \cite{SW}.

Let $\eta=(\eta_1,\eta_2,\dotsc,\eta_t)$ be a sequence of positive integers
summing to $n$ and $\gamma=(\gamma_1,\dotsc,\gamma_n)\in\ZZ^n$
an integral weight.
Given the pair $(\eta,\gamma)$ we define a generating function
$H_{\eta,\gamma}(x;q)$ for the polynomials $\K_{\la;R}(q)$.
This generating function is the graded Euler characteristic character
of the coordinate ring of a nilpotent conjugacy class closure, twisted by
a line bundle.

Let $\Roots_\eta$ be the set of positions in an $n\times n$ matrix above
the block diagonal given by the parts of $\eta$, that is,
\begin{equation*}
  \Roots_\eta = \{(i,j) \mid
	1 \le i \le \eta_1+\dotsb+\eta_r < j \le n
	\text{ for some $1\le r<t$ }\}.
\end{equation*}
Define
\begin{equation*}
  B_{\eta,\gamma}(x;q) = x^\gamma \prod_{(i,j)\in \Roots_\eta}
  	(1-q\, x_i/x_j)^{-1},
\end{equation*}
which is a formal power series in $q$ whose coefficients are
Laurent polynomials in the set of variables $x=(x_1,x_2,\dotsc,x_n)$.
Let $W$ be the symmetric group that permutes the $x$ variables,
$J=\sum_{w\in W} (-1)^w w$ the antisymmetrization operator, and
$\rho=(n-1,n-2,\dotsc,1,0)$.  Then define the formal power series
$H_{\eta,\gamma}(x;q)$ and $K_{\la,\eta,\gamma}(q)$ by
\begin{equation*}
\begin{split}
  H_{\eta,\gamma}(x;q) &= J(x^\rho)^{-1} J(x^\rho B_{\eta,\gamma}(x;q)) \\
  H_{\eta,\gamma}(x;q) &= \sum_\la \K_{\la,\eta,\gamma}(q) s_\la(x)
\end{split}  
\end{equation*}
where $\la$ runs over the dominant integral weights
and $s_\la(x) = J(x^\rho)^{-1} J(x^\rho x^\la)$ is the irreducible
character of highest weight $\la$.  It can be shown that the
formal power series $\K_{\la,\eta,\gamma}(q)$ are polynomials
with integer coefficients.  These can be calculated explicitly using
the $q$-Kostant partition formula that follows immediately
from this definition.

For the rest of the paper we shall only be concerned with the
following special case.
Let $R=(R_1,R_2,\dotsc,R_t)$ be a sequence of rectangular
partitions where $R_i$ has $\eta_i$ rows and $\mu_i$ columns.
Let $n=\sum_i \eta_i$ and
\begin{equation*}
  \gamma(R) = (\mu_1^{\eta_1},\mu_2^{\eta_2},\dotsc,\mu_t^{\eta_t}),
\end{equation*}
the weight obtained by juxtaposing the parts of the rectangular
partitions $R_1$ through $R_t$.  Say that $R$ is \textit{dominant}
if $\gamma(R)$ is, that is, the number of columns of the rectangles
in $R$ weakly decrease.  Let
\begin{equation*}
  \K_{\la;R}(q) = \K_{\la,\eta,\gamma(R)}(q).
\end{equation*}

\subsection{Notation for rectangle sequences}
\label{rect notation}
Let us fix notation associated with the sequence of rectangular
partitions $R=(R_1,R_2,\dotsc,R_t)$.  Let
$R_i$ have $\eta_i$ rows and $\mu_i$ columns for $1\le i\le t$ and let
\begin{equation*}
\begin{split}
  N &= \sum_{i=1}^t |R_i| = \sum_{i=1}^t \mu_i \eta_i \\
  n &= \sum_{i=1}^t \eta_i.
\end{split}
\end{equation*}
Let $A_1$ be the first $\eta_1$ elements of
$[n]=\{1,2,\dotsc,n\}$, $A_2$ the next $\eta_2$ elements, etc.

Let $Y_i$ be the unique column-strict tableau of shape $R_i$
and content $R_i$ in the alphabet $A_i$, meaning that the $j$-th
row of $Y_i$ consists of $\mu_i$ copies of the $j$-th largest letter
(namely $\eta_1+\eta_2+\dotsb+\eta_{i-1}+j$) of $A_i$.

\begin{ex} Let $R=((3,3),(2,2),(1,1,1))$, so that $\mu=(3,2,1)$,
$\eta=(2,2,3)$, $n=7$, $N=13$, $A_1=[1,2]$, $A_2=[3,4]$, $A_3=[5,7]$, and 
\begin{equation*}
  Y_1 = \begin{matrix} 1&1&1\\2&2&2 \end{matrix} \qquad
  Y_2 = \begin{matrix} 3&3\\4&4 \end{matrix} \qquad
  Y_3 = \begin{matrix} 5\\6\\7 \end{matrix}
\end{equation*}
\end{ex}

\begin{rem} Suppose $\eta_i=1$ for all $i$ and $\mu$ is a partition of $N$;
this shall be referred to as the Kostka case.  Then $A_i=\{i\}$
for all $i$ and $Y_i$ is the one-row tableau consisting of $\mu_i$
copies of the letter $i$.
\end{rem}

\subsection{$R$-LR words}
Let $w\in[n]^*$ be a word in the alphabet $[n]$.
For $B\subseteq[n]$, denote by $w|_B$ the restriction of $w$ to the
subalphabet $B$, obtained from $w$ by erasing all letters not in $B$.
The Knuth equivalence on words will be denoted $v\Knuth w$ \cite{Knu}.
For a (possibly skew) tableau $T$,
let $\word(T)$ denote the row-reading word of $T$, given by
$\dotsb u^2 u^1$ where $u^i$ is the $i$-th row of $T$ read from left
to right.  Often we write $T$ instead of $\word(T)$ when only the
Knuth equivalence class of $T$ is needed.

Say that $w$ is \textit{$R$-LR} (short for $R$-Littlewood-Richardson) if 
\begin{equation*}
  w|_{A_i} \Knuth Y_i \qquad\text{for all $i$}.
\end{equation*}
Denote by $W(R)$ the set of $R$-LR words.
In the Kostka case, $w$ is $R$-LR if and only if it has content $\mu$.

\begin{ex} In the running example, the word $w=7442632512131$ is $R$-LR
since $w|_{A_1}=221211\Knuth 222111$, $w|_{A_2}=4433$ and $w|_{A_3}=765$.
\end{ex}

\begin{rem} \label{rect rem} 
$W(R)$ is the set of shuffles of the Knuth classes of
$\{\word(Y_i)\}$.  Since each $\word(Y_i)$
consists of letters in the alphabet $A_i$ and the intervals $A_i$ are
non-overlapping, it follows that $W(R)$ is closed under Knuth equivalence.

Let $1\le i\le j\le t$.  Write $R_{[i,j]}=(R_i,R_{i+1},\dotsc,R_j)$
and $B := \cup_{r=i}^j A_r$.  As an immediate consequence of the
definitions, if $w\in W(R)$ then $w|_B\in W(R_{[i,j]})$.
\end{rem}

\subsection{Littlewood-Richardson tableaux}
Let $\LRT(R)$ be the set of column-strict tableaux $T$ of arbitrary
partition shape such that $\word(T)$ is $R$-LR.  For a skew shape $D$,
write $\LRT(D;R)$ for the set of column-strict tableaux $T$ of shape $D$
such that $\word(T)$ is $R$-LR.

\begin{rem} \label{LR rem} For a word $w$, denote by $P(w)$ its
Schensted $P$ tableau, the unique column-strict tableau of
partition shape such that $\word(P(w)) \Knuth w$.  By Remark \ref{rect rem},
a word $w$ is $R$-LR if and only if $P(w)\in\LRT(R)$.
\end{rem}

In the Kostka case, $\LRT(R)$ is the set of column-strict
tableaux of content $\mu$ and arbitrary partition shape, and $\LRT(D;R)$
is the set of column-strict tableaux of shape $D$ and content $\mu$.

\begin{prop} \label{rect LR coef} The multiplicity
\begin{equation*}
  \LRC^\la_R := \inner{s_\la}{s_{R_1} s_{R_2} \dotsb s_{R_t}}
\end{equation*}
is given by the cardinality of the set $\LRT(\la;R)$.
\end{prop}

This is an easy consequence of the well-known rule of Littlewood and
Richardson \cite{LR} (see section \ref{LR sec}).

\subsection{Generalization of automorphisms of conjugation}
In the Kostka case, there is a content-permuting
action of the symmetric group $S_t$ on words in the alphabet $[t]$
induced by the automorphisms of conjugation \cite{LS2} (see section
\ref{crystal}).  We recall from \cite{KS} a construction that generalizes
these bijections for $R$-LR words and LR tableaux.  It should be mentioned
that in \cite{KS} many of the properties of the generalized automorphisms
of conjugation were conjectural but are proven here using the direct
definition of cyclage on LR tableaux.

The symmetric group $S_t$ has an obvious action
on the set of all sequences of rectangles of length $t$.
For any two sequences $R$ and $R'$ in the same $S_t$-orbit,
Proposition \ref{rect LR coef} implies that
\begin{equation*}
  |\LRT(\la;R)| = |\LRT(\la;R')|.
\end{equation*}
We wish to define a family of bijections
\begin{equation*}
  \tau^{R'}_{R}:\LRT(\la;R)\rightarrow \LRT(\la;R')
\end{equation*}
that is functorial in the sense that $\tau^R_R$ is the identity
on $\LRT(\la;R)$ and
\begin{equation*}
  \tau^{R"}_R = \tau^{R"}_{R'} \circ \tau^{R'}_R
\end{equation*}
for $R$, $R'$, and $R"$ in the same $S_t$ orbit.  Let $\RRR$ be an $S_t$
orbit of a sequence of $t$ rectangles.  This functoriality then
makes it possible to define an action of the symmetric group
$S_t$ on the union $\bigcup_{R\in \RRR} \LRT(\la;R)$
such that the action of the permutation $\sigma\in S_t$
is defined by
\begin{equation*}
  \sigma|_{\LRT(\la;R)} = \tau^{\sigma R}_R
\end{equation*}
for every $R\in \RRR$.  This action of $S_t$ generalizes the
automorphisms of conjugation.

It is helpful to describe the above situation rather carefully,
since each sequence of rectangles $R\in \RRR$ has a different
set of subalphabets $\{A_1,\dotsc,A_t\}$ of $[n]$.
When the context is clear
we will abuse notation by writing $\sigma$ instead of $\tau^{\sigma R}_R$.

Recall that the (row-insertion) Robinson-Schensted (RS) correspondence 
\cite{Sch} assigns to each word $w=w_1 w_2\dotsb w_N$ (where $w_i$ is a
letter) the pair of tableaux $(P(w),Q(w))$ of the same partition shape, 
where $Q(w)$ is the unique standard tableau such that the shape of
the tableaux $Q(w)|_{[i]}$ and $P(w_1 \dotsb w_i)$ coincide for
all $0\le i\le N$.

As with the automorphisms of conjugation, the action of $S_t$
is defined on words.  By definition, the RS corespondence
restricts to a bijection
\begin{equation*}
  W(R) \rightarrow \bigcup_{\la} \LRT(\la;R) \times \ST(\la)
\end{equation*}
sending $w$ to the tableau pair $(P(w), Q(w))$, where $\ST(\la)$ is the
set of standard tableaux of shape $\la$.  We wish to define
a functorial family of bijections $\tau^{R'}_R:W(R)\rightarrow W(R')$
(for $R$ and $R'$ in $\RRR$) with the property that the following
diagram commutes:
\begin{equation*}
\begin{CD}
  W(R) @>{RS}>> \bigcup_\la \LRT(\la;R) \times \ST(\la) \\
  @V{\tau^{R'}_R}VV 	@VV{\cup_\la \tau^{R'}_R \times id_{\ST(\la)}}V \\
  W(R') @>{RS}>> \bigcup_\la \LRT(\la;R') \times \ST(\la) \\
\end{CD}
\end{equation*}
that is, $P(\tau^{R'}_R w) = \tau^{R'}_R P(w)$ and
$Q(\tau^{R'}_R w) = Q(w)$ for all $w\in W(R)$.  
Define the action of $S_t$ on $\bigcup_{R\in\RRR} W(R)$ by
\begin{equation*}
 \sigma|_{W(R)} = \tau^{\sigma R}_R.
\end{equation*}
The commutation of the above diagram can then be rephrased as
\begin{equation*}
\begin{split}
  P(\sigma w) &= \sigma P(w) \\
  Q(\sigma w) &= Q(w)
\end{split}  
\end{equation*}
for all $w\in W(R)$ and $\sigma\in S_t$.

As in the definition of the automorphisms of conjugation,
we define the action of $S_t$ in terms of the
action of the adjacent transpositions.  The simple reflections
in the symmetric group $S_t$ shall be denoted $\tau_p$ for $1\le p\le t-1$.
The corresponding bijection on LR words and tableaux
shall also be denoted $\tau_p$.  When we
have occasion to use the original automorphisms of conjugation
acting in the symmetric group $S_n$, they will be denoted
$s_r$ for $1\le r\le n-1$.  We call the bijections that give the
action of $\tau_p$ ``rectangle-switching bijections".

The definition of $\tau_p$ proceeds in a sequence of cases.  Let
\begin{equation*}
R' = \tau_p R = (R_1,\dotsc,R_{p-1},R_{p+1},R_p,R_{p+2},\dotsc,R_t)
\end{equation*}
with rectangles $(R_1',R_2',\dotsc,R_t')$,
subalphabets $A_1'$, $A_2'$, etc., and tableaux $Y_1'$, $Y_2'$, etc.

\begin{enumerate}
\item[(1)] Suppose $t=2$ and $p=1$.  Then define
$\tau_1:\LRT(\la;R)\rightarrow \LRT(\la;R')$ to be the unique map
between the two sets.  It is necessarily a bijection since both the
domain and range are either both empty or singleton sets (see Proposition
\ref{two rect}).
\item[(2)] Suppose $t=2$ and $p=1$ with the notation of the previous case.
Define the map $\tau_1:W(R)\rightarrow W(R')$ by
$P(\tau_1 w) = \tau_1 P(w)$ and $Q(\tau_1 w)=Q(w)$, which is well-defined
by the bijectivity of the RS correspondence and case (1).
\item[(3)] Let $t$ be arbitrary and $1\le p\le t-1$.
Let $B=A_p\cup A_{p+1}$.  Given $w\in W(R)$,
let $\tau_p w$ be the word obtained from $w$ by replacing the
letters of the subword $w|_B$ (in the same positions) by those of
$\tau_p (w|_B)$, which is defined by case (2), since
$w|_B\in W((R_p,R_{p+1}))$ in the alphabet $B$ by Remark \ref{rect rem}.
\item[(4)] With the same hypotheses as case (3), let $T$ be a
(possibly skew) tableau such that $\word(T)\in W(R)$.
Define $\tau_p T$ to be the unique tableau of the same shape as $T$
such that
\begin{equation*}
  \word(\tau_p T) = \tau_p \word(T).
\end{equation*}
This tableau is column-strict (Theorem \ref{action} (A1)).
\end{enumerate}

\begin{ex} With $R$ as before, consider the tableau $S\in \LRT(R)$
and some of its images under operators $\tau_p$.
\begin{equation*}
\begin{split}
  S&=\begin{matrix} 1&1&1&3&3&5\\2&2&2&4& & \\4&6& & & & \\7& & & & &
  \end{matrix} \qquad
  \tau_1 S=\begin{matrix} 1&1&3&3&3&5\\2&2&4&4& & \\4&6& & & & \\7& & & & &
  \end{matrix} \\
  \tau_2 S&=\begin{matrix} 1&1&1&3&6&6\\2&2&2&4& & \\5&7& & & & \\7& & & & &
  \end{matrix} \qquad
  \tau_2 \tau_1 S=\begin{matrix} 1&1&3&6&6&6\\2&2&4&7& & \\5&7& & & & \\7& & & & &
  \end{matrix} \\
  \tau_1 \tau_2 S&=\begin{matrix} 1&4&4&4&6&6\\2&5&5&5& & \\3&7& & & & \\7& & & & &
  \end{matrix} \qquad
  \tau_1 \tau_2 \tau_1 S=
	\begin{matrix} 1&4&4&6&6&6\\2&5&5&7& & \\3&7& & & & \\7& & & & &
  \end{matrix}
\end{split}
\end{equation*}
\end{ex}

\begin{rem} \label{auto}
In the Kostka case the bijection $\tau_p$ is precisely the 
automorphism of conjugation $s_p$.
\end{rem}

Say that a word $w$ \textit{fits} the skew shape $D$ if there
is a column-strict tableau $T$ of shape $D$ such that $w=\word(T)$.

\begin{thm} \label{action} The maps $\tau_p$ are well-defined bijections
\begin{equation*}
\begin{split}
W(R)&\rightarrow W(\tau_p R) \\
\LRT(\la;R)&\rightarrow\LRT(\la;\tau_p R)
\end{split}
\end{equation*}
satisfying the following properties.
\begin{enumerate}
\item[(A1)] If $D$ is a skew shape, then $w\in W(R)$ fits $D$ if and only
if $\tau_p w$ fits $D$.
\item[(A2)] If $v\Knuth w$ are in $W(R)$, then
$\tau_p v\Knuth \tau_p w$ in $W(\tau_p R)$.
\item[(A3)] $P(\tau_p w) = \tau_p P(w)$ for all $w\in W(R)$.
\item[(A4)] $Q(\tau_p w) = Q(w)$ for all $w\in W(R)$.
\item[(A5)] The bijections $\tau_p$ ($1\le p\le t-1$) satisfy the
Moore-Coxeter relations for $S_t$, defining an action of the symmetric
group on LR tableaux and words.
\item[(A6)] Suppose $\sigma\in S_t$ stabilizes the subintervals
$[1,i-1]$, $I=[i,j]$, and $[j+1,t]$.  Let $w\in W(R)$ and let
$\sigma w$ denote the action of $\sigma \in S_t$ on $w$ given in (A5).
Then the positions of the letters of $J=\bigcup_{i\in I} A_j$ in $w$
coincide with the positions of the letters of $J$ in $\sigma w$.
\item[(A7)] Suppose $R_p=R_{p+1}$.  Then $\tau_p(w)=w$ for all
$w\in W(R)$.
\end{enumerate}
\end{thm}

\subsection{Statistic on LR tableaux}
Lascoux, Leclerc, and Thibon gave a formula for the charge
that averages a simpler statistic over the orbit of a word
by the symmetric group acting by automorphisms of conjugation
\cite{LLT2}.  We recall from \cite{KS} a statistic $\LLT_R$
on LR tableaux that generalizes the charge.  In \cite{KS}
even the well-definedness of $\LLT_R$ was based on a conjecture
(namely, Theorem \ref{action} (A5)); here these conjectures
regarding $\LLT_R$ are proven.

First let $R=(R_1,R_2)$ and $w\in W(R)$.  Then $P(w)\in\LRT(R)$.
Define
\begin{equation*}
  d_{R_1,R_2}(w) = d_{R_1,R_2}(P(w))
\end{equation*}
to be the number of cells of the shape of $P(w)$ that lie in columns
strictly to the right of the $\max(\mu_1,\mu_2)$-th column.
In the notation of Proposition \ref{two rect}, this is the size of
the partition $\la_{ne}$ where $\la$ is the shape of $P(w)$.

Next, for arbitrary $R$, $1\le i\le t-1$, and $w\in W(R)$, recall that
by definition $w|_{A_i\cup A_{i+1}}\in W((R_i,R_{i+1}))$.
Define $d_{i,R}(w) = d_{(R_i,R_{i+1})}(w|_{A_i\cup A_{i+1}})$.

Finally, define the statistic $\LLT_R:W(R)\rightarrow\NN$ as
the following average over the symmetric group $S_t$:
\begin{equation} \label{LLT charge}
  \LLT_R(w) = \dfrac{1}{t!} \sum_{\tau\in S_t}
  	\sum_{i=1}^{t-1}\, (t-i)\, d_{i,\tau R}(\tau w).
\end{equation}
This makes sense if $R$ has two or more rectangles.
Let $\LLT_R(w)=0$ if $R$ has less than two rectangles.
For $T\in\LRT(R)$, define $\LLT_R(T)=\LLT_R(\word(T))$.

\begin{ex} For the tableaux in the $S_3$-orbit of the tableau $S$
in the previous example, in order one has the following
values for $d_1$ and $d_2$: $(3,1)$, $(3,1)$, $(2,1)$, $(2,1)$,
$(2,2)$, $(2,2)$.  So $\LLT_R(T)=1/6 (7+7+5+5+6+6) = 6$.
\end{ex}

\subsection{Main result}

\begin{thm} \label{main result} Let $R$ be a dominant sequence of
rectangles and $\la$ a partition.  Then
\begin{equation} \label{LR poly}
  \K_{\la;R}(q) = \sum_{T\in\LRT(\la;R)} q^{\LLT_R(T)}.
\end{equation}
\end{thm}

The proof relies on the cyclage poset structure on LR tableaux
introduced below.

\section{The cyclage poset for LR tableaux}
\label{cyclage sec}

Section \ref{crystal} reviews the definitions of the crystal
operators for type $A_{n-1}$, including the action of the
symmetric group $S_n$ given by the automorphisms of conjugation.
This is necessary to define the action of the cyclic group
$\ZZ/N\ZZ$ on $R$-LR words in section \ref{cyclic action}, 
which in turn makes possible the definition in section \ref{cocyclage sec}
of the $R$-cocyclage relation on $R$-LR tableaux.  Section
\ref{cocyclage sec} gives the main theorems on the structure
of the $R$-cocyclage poset on $R$-LR tableaux and its grading function
$\charge_R$.

\subsection{Crystal operators and automorphisms of conjugation}
\label{crystal}

We recall the definitions of the crystal raising, lowering, and
reflection operators $e_r$, $f_r$, and $s_r$ associated with the crystal
graph of type $A_{n-1}$.  The combinatorial constructions as given here
appear in \cite{LS2} (see also \cite{LLT2}),
and are equivalent to those coming from the computation
of this crystal graph by Kashiwara and Nakashima \cite{KN}.

Fix $1\le r\le n-1$.  For a word $u$, regard the letters $r$ as
right parentheses, the letters $r+1$ as left parentheses, and ignore
other letters.  Perform the usual matching of parentheses, leaving
a subword of unpaired letters of the form $r^a (r+1)^b$.  Then the words
$e_r u$, $f_r u$, and $s_r u$ are defined by replacing this subword of $u$
(in the same positions) by $r^{a+1} (r+1)^{b-1}$,
$r^{a-1} (r+1)^{b+1}$, and $r^b (r+1)^a$ respectively, where
$e_r u$ and $f_r u$ are defined only if $b>0$ and $a>0$ respectively.
Say that two words are in the same $r$-string if one is obtained from
the other by a power of $f_r$.

\begin{thm} \cite{LS2} There is an action of $S_n$ on the words in the
alphabet $[n]$, where the simple reflection $s_r\in S_n$ acts
by the above operator $s_r$ on words.
\end{thm}

The operators on words corresponding to permutations of $S_n$
are called automorphisms of conjugation in \cite{LS2}.

\subsection{Action of the cyclic group $\ZZ/N\ZZ$ on $W(R)$}
\label{cyclic action}

In the Kostka case, the cyclic group $\ZZ/N\ZZ$ acts on words of content
$\mu$ by cyclic rotation of positions, leading to the definition of the
cyclage poset on column-strict tableaux of content $\mu$ \cite{LS2}. 
This simple action of $\ZZ/N\ZZ$ is extended to the set of $R$-LR words.

\begin{ex} The naive action that merely rotates positions does not 
preserve the set of $R$-LR words.  It is illustrative to
consider the case $R=(R_1)$.  Write $a=\mu_1$ and $Y=Y_1$.
Consider $w=\word(Y)$ and its right circular rotation $v$,
given explicitly by
\begin{equation*}
\begin{split}
  w &= n^a (n-1)^a \dotsb 2^a 1^a \\
  v &= 1\, n^a (n-1)^a \dotsb 2^a 1^{a-1} \\
\end{split}
\end{equation*}
Clearly $v$ is not $R$-LR unless $n=1$.
\end{ex}

This right circular rotation must be modified to preserve $R$-LRness.
Let $w\in[n]^*$.  Write $w=ux$ where $x$ is a letter.  Let
$w_0^R$ be the longest permutation in the Young subgroup of
$S_n$ that stabilizes all the intervals $A_i$.  Define
\begin{equation}
  \chi_R(w) = (w_0^R x)(w_0^R u)
\end{equation}
where $w_0^R$ acts by the automorphism of conjugation
(see Remark \ref{auto}).

In the Kostka case, the interval $A_i=\{i\}$ so $w_0^R$ is the
identity permutation and $\chi_R(w)$ is merely the right circular
rotation of $w$.

\begin{ex} In the previous example, $w_0^R$ is the longest permutation
in $S_n$, $x=1$, $x w_0^R = n$,
\begin{equation*}
\begin{split}
  u &= n^a(n-1)^a\dotsb 2^a 1^{a-1} \\
  u w_0^R &= (n-1) n^{a-1} (n-2) (n-1)^{a-1} \dotsb 2 3^{a-1}
		1 2^{a-1} 1^{a-1}.
\end{split}
\end{equation*}
Observe that $\chi_R(w)\in W(R)$.
\end{ex}

\begin{prop} \label{chi} Let $k$ be a nonnegative integer.
\begin{enumerate}
\item $w\in W(R)$ if and only if $\chi_R^k(w)\in W(R)$.
\item Let $w=uv\in W(R)$ with $v$ of length $k$.
respectively.  Then
\begin{equation*}
  \chi_R^k(uv)=(w_0^R v)(w_0^R u)
\end{equation*}
where $w_0^R\in S_n$ acts by an automorphism of conjugation.
\end{enumerate}
\end{prop}

The following key result generalizes \cite[Theorem 4.6]{LS2}.

\begin{thm} \label{cocyc action} Let $R'=\tau_p R$.  The following
diagram commutes:
\begin{equation*}
\begin{CD}
  W(R) @>{\chi_R}>> W(R) \\
  @V{\tau_p}VV		@VV{\tau_p}V \\
  W(R')@>>{\chi_{R'}}> W(R')
\end{CD}
\end{equation*}
\end{thm}

\subsection{$R$-cocyclage}
\label{cocyclage sec}

Let $S,T\in \LRT(R)$.  Let $\le_R$ be the transitive closure
of the relation $T<_R S$, which holds if there is a word $w=ux\in W(R)$
with $x$ a letter, such that $S=P(w)$, $T=P(\chi_R(w))$,
and the cell given by the difference of the shapes of $P(w)$
and $P(u)$, is in a column strictly east of the $a$-th
where $a=\max_i \mu_i$.

\begin{rem} \label{corner chi}
All covering relations $T<_R S$ can be realized as follows.  Let
$S\in \LRT(R)$ and $s$ a corner cell of $S$ in a column
strictly east of the $a$-th.  Perform the reverse row insertion
on $S$ starting at $s$, producing the tableau $U$ and ejecting the
letter $x$.  Then letting $u=\word(U)$ and $w=ux$, one has
$T<_R S$ where $T=P(\chi_R(w))$.
\end{rem}

\begin{ex} Using the tableau $S$ from the running example
and the corner cell $(2,4)$, we have $w_0^R = s_1 s_3 s_5 s_6 s_5$.
\begin{equation*}
S=\begin{matrix} 1&1&1&3&3&5\\2&2&2&4& & \\4&6& & & &
		\\7& & & & & \end{matrix} \qquad
U=\begin{matrix} 1&1&1&3&4&5\\2&2&2& & & \\4&6& & & &
		\\7& & & & & \end{matrix}
\end{equation*}
and $x=3$.  Now $w_0^R x = 4$ and 
\begin{equation*}
w_0^R U = \begin{matrix} 1&1&1&3&3&5\\2&2&2& & & \\4&6& & & &
		\\7& & & & & \end{matrix} \qquad
T=\begin{matrix} 1&1&1&3&3&5\\2&2&2& & & \\4&4&6& & &
		\\7& & & & & \end{matrix}
\end{equation*}
\end{ex}

\begin{thm} \label{poset thm}
\begin{enumerate}
\item $(\LRT(R),\le_R)$ is a graded poset.
\item An element of $\LRT(R)$ is $\le_R$-minimal if and only if it has
exactly $a=\max_i \mu_i$ columns.  In particular, if
all of the rectangles in $R$ have the same number of columns then
then there is a unique minimal tableau.
\item Suppose $\mu_1=a$ and $T\in\LRT(R)$ is $\le_R$-minimal.  Then
$T$ consists of the tableau $Y_1$ sitting atop a tableau
$\Th\in\LRT(\Rhat)$ in the alphabet $[\eta_1+1,n]$ where
$\Rhat=(R_2,R_3,\dotsc,R_t)$.
\end{enumerate}
\end{thm}

\begin{thm} \label{switch iso} $\tau_p$ is an isomorphism of
graded posets
\begin{equation*}
  (\LRT(R),\le_R)\rightarrow (\LRT(\tau_p R),\le_{\tau_p R})
\end{equation*}
\end{thm}

\begin{thm} \label{intrinsic charge_R} There is a unique function
$\charge_R:\LRT(R)\rightarrow\NN$ such that:
\begin{enumerate}
\item[(C1)] If $R=\emptyset$, then $\charge_R(\emptyset)=0$.
\item[(C2)] If $T' <_R T$ is a covering relation
for $T$ and $T'$ in $\LRT(R)$, then $\charge_R(T)=\charge_R(T')+1$.
\item[(C3)] With the assumptions and notation of Theorem \ref{poset thm} 3,
\begin{equation*}
	\charge_R(T)=\charge_{\Rhat}(\Th).
\end{equation*}
\item[(C4)] $\charge_{\tau_p R}(\tau_p T) = \charge_R(T)$ for $T\in\LRT(R)$.
\end{enumerate}
\end{thm}

\begin{thm} \label{explicit charge} The function $\LLT_R$ is constant on
Knuth equivalence classes and satisfies conditions
(C1) through (C4) of Theorem \ref{intrinsic charge_R}.
\end{thm}

\begin{prop} \label{chi charge} Let $w\in W(R)$.  For each $\tau\in S_t$,
let $f(\tau)\in[t]$ be such that the last letter of $\tau w$
is in the alphabet $A_{f(\tau)}$, and $\NC$ the set of indices $i\in[t]$
such that $f(\tau_1 \tau_2 \dotsb \tau_{i-1})=1$.  Then
\begin{equation} \label{LLT diff}
  \LLT_R(w) = \LLT_R(\chi_R(w)) + 1 - |\NC|.
\end{equation}
\end{prop}

\section{Proof of main theorem}
\label{proof outline}

The proof follows the suggested outline in \cite{SW}.
The important difference is that here we work with
LR tableaux, as opposed to the $R$-catabolizable
tableaux used in \cite{SW}.  As in the approach of
Lascoux and Sch\"utzenberger for the Kostka-Foulkes polynomials,
Theorem \ref{main result} is proven by showing that the right 
hand side of \eqref{LR poly} satisfies (a form of) the generalized Morris
recurrence for $\K_{\la;R}(q)$ \cite{SW}:
\begin{equation}
\label{new rec}
  \K_{\la;R}(q) = \sum_{w\in W} (-1)^w q^{|\alpha(w)|-|R_1|}
  	\sum_\nu K_{\nu,(\alpha(w)-R_1,\beta(w))}
  		\K_{\nu;\Rhat}(q)
\end{equation}
where $\alpha(w)$ and $\beta(w)$ are the first $m:=\eta_1$ and last
$n-m$ parts of the weight
\begin{equation*}
  \xi(w) := w^{-1}(\la+\rho)-\rho
\end{equation*}
and $K_{\la,\alpha}$ is the Kostka number \cite[I.6]{Mac}, the
multiplicity of the weight $\alpha$ in the $\la$-th irreducible
$gl(n)$-module, or equivalently the number of column strict tableaux
of shape $\la$ and content $\alpha$.

Let $\SSS$ be the set of triples $(w,T,U)$ where $w\in W$, and
$T$ and $U$ are column-strict tableaux of the same partition shape,
where the word of $T$ is $\Rhat$-LR in the alphabet $[m+1,n]$
and $U$ has content $(\beta(w),\alpha(w)-R_1)$. The reordering of parts of
the content of $U$ is justified since the Kostka number is symmetric in
its second index.  Define a sign and weight on $\SSS$ by
\begin{equation*}
\begin{split}
  \sign(w,T,U) &= (-1)^w \\
  \weight(w,T,U) &= q^{|\alpha(w)|-|R_1| + \charge_{\Rhat}(T)}.
\end{split}
\end{equation*}
By induction the right hand side of \eqref{new rec} is given by
\begin{equation} \label{sign weight}
  \sum_{(w,T,U)\in\SSS} \sign(w,T,U) \weight(w,T,U).
\end{equation}

We wish to write the sum \eqref{sign weight} as the generating function
of another signed weighted set that is more amenable to cancellation.
Let $\TTT$ be the set of triples $(w,P,Q)$ such that $w\in W$
and $P$ and $Q$ are column-strict tableaux of the same partition shape,
such that the word of $P$ is $R$-LR and $Q$ has content $\xi(w)$.
Define
\begin{equation*}
\begin{split}
  \sign(w,P,Q) &= (-1)^w \\
  \weight(w,P,Q) &= q^{\charge_R(P)}.
\end{split}
\end{equation*}
We define a sign- and weight-preserving map
$\Phi:\SSS\rightarrow\TTT$ as follows.

The column insertion version of the Robinson-Schensted-Knuth correspondence
(column RSK) sends a sequence of weakly increasing words
$(u^n,u^{n-1},\dotsc,u^1)$ to a pair of column-strict tableaux $(P,Q)$ where
$P=P(u^n u^{n-1}\dotsb u^1)$ is the usual Schensted $P$-tableau 
and $Q$ is the column-strict tableau in the alphabet $[n]$ such that
the shape of $Q|_{[i]}$ coincides with the shape of
$P(u^i u^{i-1} \dotsb u^1)$ for all $0\le i\le n$.

Given the triple $(w,T,U)\in\SSS$, consider the sequence of words
given by the inverse image of the tableau pair $(T,U)$ under column RSK.
Here we use a nonstandard indexing for the weakly increasing words, writing
\begin{equation*}
P(u^m u^{m-1} \dotsb u^1 u^n u^{n-1} \dotsb u^{m+1}) = T.
\end{equation*}
Since $U$ has content $(\beta(w),\alpha(w)-R_1)$,
this nonstandard indexing allows us to say that the length of
$u^i$ is the $i$-th part of the weight
$\xi(w)-(R_1,0^{n-m}) = (\alpha(w)-R_1,\beta(w))$.

Let $y^i$ be the word of the same length as $u^i$, such that
\begin{equation*}
  \chi_{\Rhat}^{|\alpha(w)|-|R_1|}(y^n y^{n-1} \dotsb y^1) =
	u^m u^{m-1} \dotsb u^1 u^n u^{n-1} \dotsb u^{m+1}.
\end{equation*}
This is well-defined since $\chi_{\Rhat}$ is invertible.
By Proposition \ref{chi} and Theorem \ref{action} for
the automorphism of conjugation $w_0^{\Rhat}$, it follows that
each $y^i$ is a weakly increasing word with letters greater than $m$.  Let
\begin{equation} \label{v def}
  v^i = \begin{cases}
  	i^{\mu_1} y^i & \text{if $1\le i\le m$} \\
  	y^i	& \text{if $m<i\le n$.}
  	\end{cases}
\end{equation}
The word $v^i$ is weakly increasing since $y^i$ is.
Note that $|v^i| = \xi(w)_i$.

Finally, let $(P,Q)$ be the pair of tableaux given by the image under
column RSK, of the sequence of words $\{v^i\}$ so that
$P(v^n v^{n-1} \dotsb v^1) = P$.  Note that $Q$ has content $\xi(w)$.
Then define $\Phi(w,T,U):=(w,P,Q)$.

\begin{ex} Let $n=8$, $\eta=(2,2,2,1,1)$, and $\mu=(3,3,3,3,3)$, so that
$m=2$ and
\begin{equation*}
\begin{split}
  R&=((3,3),(3,3),(3,3),(3),(3)) \\
  \Rhat&=((3,3),(3,3),(3),(3)).
\end{split}
\end{equation*}
Let $w=32154678$ and $\la=(6,5,5,5,2,1,0,0)$\, so that
$\xi(w) = (3,5,8,1,6,1,0,0)$, $\alpha(w)=(3,5)$ and
$\beta(w)=(8,1,6,1,0,0)$.  Let $T$ and $U$ be given by
\begin{equation*}
  T = \begin{matrix}
    3&3&3&5&5&7&7&7\\
    4&4&4&6&8&8& & \\
    5&6&8& & & & & \\
    6& & & & & & & 
  \end{matrix} \qquad
  U = \begin{matrix}
    1&1&1&1&1&1&1&1\\
    2&3&3&3&3&3& & \\
    3&8&8& & & & & \\
    4& & & & & & & 
  \end{matrix}
\end{equation*}
The words $u^i$ are given by
$u^2=66$, $u^1=u^8=u^7=\emptyset$, $u^6=8$,
$u^5=445688$, $u^4=4$, $u^3=33355777$.
Since
\begin{equation*}
  \chi_{\Rhat}^{-2}(66/ / / /8/445688/4/33355777) =
  	   / /8/446688/4/33356777/55/,
\end{equation*}
we have $v^1=111$, $v^2=22255$,
$v^3=33356777$, $v^4=4$, $v^5=446688$, $v^6=8$, and $v^7=v^8=\emptyset$.

The tableaux $P$ and $Q$ are given by
\begin{equation*}
P=\begin{matrix}
  1&1&1&5&5&5&7&7\\
  2&2&2&6&6&7& & \\
  3&3&3&8&8& & & \\
  4&4&4& & & & & \\
  6& & & & & & & \\
  8& & & & & & &
\end{matrix}  \qquad
Q=\begin{matrix}
1&1&1&2&2&3&3&3 \\
2&2&2&3&3&5& &  \\
3&3&3&5&5& & &  \\
4&5&5& & & & &  \\
5& & & & & & &  \\
6& & & & & & &
\end{matrix}
\end{equation*}
\end{ex}

\begin{lem} \label{preserving map} The map $\Phi:\SSS\rightarrow\TTT$
is a sign- and weight-preserving injection.
\end{lem}

By Lemma \ref{preserving map} and \eqref{sign weight}, the
right hand side of \eqref{new rec} is equal to
\begin{equation} \label{start cancel}
  \sum_{(w,P,Q)\in\Phi(\SSS)} \sign(w,P,Q) \weight(w,P,Q).
\end{equation}

We define a sign-reversing, weight-preserving involution
$\theta$ on the set $\TTT$.  Let $(w,P,Q)\in\TTT$.
\begin{enumerate}
\item If $Q$ is lattice (see section \ref{LR sec}), it follows that
$Q$ has partition content, $w$ is the identity, and the content
of $Q$ is $\la$.  In this case define $\theta(w,P,Q)=(w,P,Q)$.
\item If $Q$ is not lattice, let $r+1$ be the rightmost
letter in $\word(Q)$ that violates the lattice
property.  Define $\theta(w,P,Q)=(w',P',Q')$ where
$w' = w s_r$, $P'=P$, and $Q' = s_r e_r Q$.
\end{enumerate}

$\theta$ is an involution on all of $\TTT$, due to the
following easy lemma.

\begin{lem} \label{Bott involution}  \cite{Sh}
Consider the following map $\Psi$ on the set of non-lattice
words.  Suppose the word $u$ is not lattice.  Let $r+1$
be the rightmost letter where latticeness fails, and let
$\Psi(u)= s_r e_r u$.  Then $\Psi$ is an involution.
\end{lem}

\begin{ex} In computing $\theta'(w,P,Q)=(w',P',Q')$, the first
violation of latticeness in the row-reading word of
$Q$ occurs at the cell $(1,8)$ so $r=2$.  We have
$w'=31254678$, $P'=P$ and
\begin{equation*}
Q'=\begin{matrix}
1&1&1&2&2&2&2&3 \\
2&2&2&3&3&5& &  \\
3&3&3&5&5& & &  \\
4&5&5& & & & &  \\
5& & & & & & &  \\
6& & & & & & &
\end{matrix}
\end{equation*}
\end{ex}

\begin{lem} \label{inv lem}
$\theta$ is a sign-reversing, weight-preserving involution
on the subset $\Phi(\SSS)$ of $\TTT$.
\end{lem}

By definition the fixed points of the restriction of $\theta$
to $\Phi(\SSS)$ are the triples
$(w,P,Q)$ where $w$ is the identity, $Q=\key(\la)$ (see section \ref{LR sec}),
and $P\in\LRT(\la;R)$.  Cancelling down from \eqref{start cancel}
using $\theta$, we have
\begin{equation*}
\begin{split}
  \K_{\la;R}(q) &= \text{ right hand side of \eqref{new rec} } \\
  &= \text{ right hand side of \eqref{start cancel} } \\
  &= \sum_{P\in\LRT(\la;R)} q^{\charge_R(P)},
\end{split}
\end{equation*}
proving Theorem \ref{main result}.

\section{Proofs}

\subsection{$R$-LR property and proof of Proposition \ref{rect LR coef}}
\label{LR sec}

Say that a word is \textit{lattice} if the content of every final subword
is a partition.

More generally, if $\mu$ is a partition, say that 
a word is \textit{$\mu$-lattice} if, for every final subword,
the sum of $\mu$ and its content, is a partition.  For the skew shape
$\la/\mu$ and the partitions $\alpha$ and $\beta$, let
$\LRT(\la/\mu;\alpha,\beta)$ be the set of column-strict tableaux of shape
$\la/\mu$ and content $\alpha-\beta$ whose row-reading words are
$\beta$-lattice.

For every composition $\alpha$ (sequence of nonnegative integers,
almost zero) let $\alpha^+$ be the partition obtained by sorting
$\alpha$ into weakly decreasing order.  There is a unique
column-strict tableau of shape $\alpha^+$ and content $\alpha$,
called the \textit{key tableau} of content $\alpha$, denoted
$\key(\alpha)$.  Explicitly, the $j$-th column of $\key(\alpha)$ consists of the
letters $i$ such that $\alpha_i\ge j$.  When $\alpha$ is a partition
$\key(\alpha)$ is called a \textit{Yamanouchi} tableau.

The following result has a straightforward proof.

\begin{prop} \label{lattice and pairs}
\begin{enumerate}
\item A word is $\mu$-lattice if and only if the number of $r$-unpaired
letters $r+1$ is at most $\mu_r-\mu_{r+1}$ for all $r$.
\item The property of $\mu$-latticeness is preserved under Knuth equivalence.
\item A word $u$ is lattice of content $\la$ if and only if $\la$
is a partition and $u \Knuth \key(\la)$.
\end{enumerate}
\end{prop}

Write $\LRT(D;\alpha;\beta)$ for the set of column-strict tableaux
of shape $D$, content $\alpha$, whose row-reading words are $\beta$-lattice.

\begin{thm} \cite{LR} (Littlewood-Richardson rule) The coefficient
\begin{equation*}
  \LRC^\la_{\mu\nu} := \inner{s_\la}{s_\mu s_\nu}
\end{equation*}
is equal to the cardinality of the set $\LRT(\la/\mu;\nu;())$.
\end{thm}

The LR rule has an immediate corollary that gives the
skew-skew LR coefficient.

Given the skew shapes $D$ and $E$, let $D\otimes E$
denote any skew shape given by the union of a translate of $D$
and a translate of $E$ such that every cell of $D$ is strictly
south and strictly west of every cell of $E$.  Clearly
$s_{D\otimes E} = s_D s_E$ for any skew shapes $D$ and $E$.

\begin{cor} \label{skew skew} The coefficient
\begin{equation*}
  \LRC^{\la/\mu}_{\alpha/\beta} := \inner{s_{\la/\mu}}{s_{\alpha/\beta}}
\end{equation*}
is equal to the cardinality of the set $\LRT(\la/\mu;\alpha,\beta)$.
\end{cor}
\begin{proof} Applying the LR rule to the skew shape
$\la/\mu \otimes \beta$ and the partition $\alpha$, one has
\begin{equation*}
\begin{split}
  \inner{s_{\la/\mu}}{s_{\alpha/\beta}}
  &= \inner{s_{\la/\mu} s_\beta}{s_\alpha} \\
  &= \inner{s_{\la/\mu} \otimes \beta}{s_\alpha} \\
  &= |\LRT(\la/\mu\otimes \beta;\alpha,())|
\end{split}
\end{equation*}
Let $T'\in \LRT(\la/\mu\otimes\beta;\alpha,())$.  Write
$T' = T \otimes U$ where $T$ and $U$ are column-strict tableaux
of shape $\la/\mu$ and $\beta$ respectively.
Clearly $\word(T')=\word(T)\word(U)$.
By the definition of latticeness, any final subword of $\word(T')$ must
itself be lattice.  In particular $\word(U)$ is lattice; it is also
the row-reading word of the column-strict tableau $U$ of shape $\beta$.
It follows from Proposition \ref{lattice and pairs} that
$U=\key(\beta)$.  This given, a direct translation of the definitions 
shows that there is a bijection
\begin{equation*}
\LRT(\la/\mu;\alpha,\beta) \rightarrow
\LRT(\la/\mu\otimes \beta; \alpha,())
\end{equation*}
given by $T\mapsto T\otimes \key(\beta)$.
\end{proof}

Proof of Proposition \ref{rect LR coef}:
\begin{proof} The skew shape
$R_t \otimes R_{t-1} \otimes \dotsb\otimes R_1$ can be explicitly
written in the form $\alpha/\beta$, where
\begin{equation*}
\begin{split}
  \alpha&=((\mu_1+\dotsb+\mu_t)^{\eta_1},(\mu_2+\dotsb+\mu_t)^{\eta_2},\dotsc,
	\mu_t^{\eta_t})\\
  \beta&=((\mu_2+\dots+\mu_t)^{\eta_1},(\mu_3+\dots+\mu_t)^{\eta_2},\dotsc,
  	\mu_t^{\eta_{t-1}},0^{\eta_t}).
\end{split}
\end{equation*}
Note that the row indices corresponding to the translate of $R_i$
are given by the subinterval $A_i$ defined in section \ref{rect notation}.
By Corollary \ref{skew skew},
\begin{equation*}
  \LRC^\la_R = \inner{s_\la}{s_{R_1} s_{R_2} \dotsb s_{R_t}} =
  \inner{s_\la}{s_{R_t\otimes R_{t-1}\otimes \dotsb \otimes R_1}} =
  \inner{s_\la}{s_{\alpha/\beta}} =
  |\LRT(\la;\alpha,\beta)|.
\end{equation*}
The $\beta$-lattice condition is vacuous for indices
$r$ of the form $\eta_1+\eta_2+\dotsb+\eta_i$, since for such $r$,
$\beta_r-\beta_{r+1}=\mu_{r+1}$, which is the total number of
letters $r+1$ in a word of content $\alpha-\beta$.
For $r$ not of this form, say
$\eta_1+\dotsb+\eta_{i-1}+1 \le r<\eta_1+\dotsb+\eta_i$, we have
$\beta_r-\beta_{r+1}=0$.  From Proposition \ref{lattice and pairs}
and the definitions, by restricting to the subalphabets $A_i$,
a word is $\beta$-lattice of content $\alpha-\beta$ if
and only if it is $R$-LR.
\end{proof}

\subsection{The one-rectangle case}
For later proofs it is necessary to conduct
detailed analyses of the cases where $R$ consists of one or
two rectangles.

Suppose $R=(R_1)$.  Then $n=\eta_1$.
Write $a=\mu_1$ and $Y=Y_1=\key((a^n))$.

The following result is easy to show by direct computation
and Proposition \ref{lattice and pairs}.

\begin{prop} \label{one rect} Suppose $u$ and $v$ words
of content $\alpha$ and $\beta$ respectively. The following are equivalent.
\begin{enumerate}
\item $vu \Knuth Y$.
\item $\alpha$ is a partition, $\alpha+\beta=(a^n)$,
$u \Knuth \key(\alpha)$ and $v\Knuth \key(\beta)$.
\end{enumerate}
\end{prop}

\subsection{The two rectangle case}
An important feature of this case is that the set of tableaux
$\LRT(\la;(R_1,R_2))$ is either empty or a singleton.  If it is
the latter it can
be described explicitly.  In the Kostka case these are 
column-strict tableaux of partition shape with
letters in the alphabet $\{1,2\}$.

The pair of rectangles $R=(R_1,R_2)$ defines three disjoint regions in
$\NN^2$: the partition diagram $R_1\cup R_2$, the northeast region $\NE(R)$
consisting of the cells $(r,c)$ such that
$r \le l' := \min(\eta_1,\eta_2)$ and $c > a:=\max(\mu_1,\mu_2)$,
and the southwest region $\SW(R)$, consisting of the cells
$(r,c)$ such that $r>l:=\max(\eta_1,\eta_2)$ and
$c\le a':=\min(\mu_1,\mu_2)$.  Clearly for the reversed sequence
$R'=(R_2,R_1)$ one has $\NE(R')=\NE(R)$ and $\SW(R')=\SW(R)$.

\begin{prop} \label{two rect} The set $\LRT(\la;(R_1,R_2))$ is either
empty or a singleton.  The following are necessary and
sufficient conditions that the set $\LRT(\la;(R_1,R_2))$ be nonempty.
\begin{enumerate}
\item $R_1\cup R_2 \subset \la \subset R_1\cup R_2\cup\NE(R)\cup\SW(R)$.
\item The shapes $\la_{ne}:=\la\cap\NE(R)$ and $\la_{sw}:=\la\cap\SW(R)$
are complementary inside the rectangle $R_1\cap R_2$ in the sense that
both are contained in $R_1\cap R_2$, and if $R_1\cap R_2$ is skewed
by the removal of one, then the resulting shape is the
180 degree rotation of the other.
\end{enumerate}
\end{prop}
\begin{proof} Suppose $\LRT(\la;(R_1,R_2))$ is nonempty; let $T$ be
a member.  Let $R'=(R_2,R_1)$, with alphabets $A_1'$ and $A_2'$
and Yamanouchi tableaux $Y_1'$ and $Y_2'$.  Due to the symmetry of tensor
product multiplicities and Proposition \ref{rect LR coef},
$\LRT(\la;(R_2,R_1))$ is also nonempty, containing the tableau $T'$ say.
By definition $T|_{A_1}=Y_1$ which has shape $R_1$, so $R_1\subset\la$.
A similar argument applied to $T'$ shows that $R_2\subset\la$.
By switching the rectangles $R_1$ and $R_2$ if necessary, we may assume
that either $\eta_1>\eta_2$, or $\eta_1=\eta_2$ and $\mu_1\ge \mu_2$.  In particular
$R_1$ is not properly contained in $R_2$.

Suppose $\la\not\subset R_1\cup R_2\cup\NE(R)\cup\SW(R)$.
We consider cases.  Suppose first that $R_2$ is not contained in $R_1$,
that is, $\eta_1>\eta_2$ and $\mu_1<\mu_2$.  By assumption
$(\eta_2+1,\mu_1+1)\in\la$.  This cell and those due north of it,
are in $\la$ but not in $R_1$.  Thus the skew tableau $T-Y_1$
contains a column of length $\eta_2+1$, which means it contains at
least that many distinct letters, which is impossible since it
only contains the letters of $A_2$, which has cardinality $\eta_2$.
Otherwise suppose that $R_2$ is contained in $R_1$.
Then by assumption $\la$ either contains the cell
$(l+1,a'+1)$ or the cell $(l'+1,a+1)$.  If $(l+1,a'+1)\in\la$
then one arrives at a contradiction in a manner similar to the
previous case.  Suppose $(l'+1,a+1)\in\la$.  It and all the cells
due west of it, are contained in $\la$ but are not in $R_2$.
This means that the skew tableau $T'-Y_1'$ contains a row
of length $a+1$.  But then $\word(T'-Y_1')$ contains a
contiguous weakly increasing subword of length $a+1$,
so by Theorem \ref{Pieri}, the recording tableau
of $\word(T'-Y_1')$ contains cells in at least $a+1$ different columns.
It follows that the shape of $P(T'-Y_1')=Y_2'$ has at least
$a+1$ columns, which is a contradiction, since $Y_2'$
has $a=\mu_1$ columns.  Thus $\la$ is contained in
$R_1\cup R_2\cup\NE(R)\cup\SW(R)$.

To show the second condition is necessary, recall that the
column-reading word of a (possibly skew) tableau $S$ is
the word $u^1 u^2\dotsb$ where $u^i$ is the word comprising the $i$-th
column of $S$, read from bottom to top.  Let us consider the
column-reading word $\cword(T-Y_1)$.  By definition it satisfies
$\cword(T-Y_1)\Knuth Y_2$.  Consider another dissection of
$\la$ into three pieces: $R_1$, the part $\la_w$ of $\la-R_1$
in the first $\mu_1$ columns, and the part $\la_e$ of $\la-R_1$
not in the first $\mu_1$ columns.  Label the corresponding parts of
the tableau $T$ by $Y_1$, $T_w$, and $T_e$.
Clearly $\cword(T-Y_1)=\cword(T_w)\cword(T_e)$.
By Proposition \ref{one rect}, $T_w=\key(\beta)$ and
$T_e=\key(\alpha)$ (both in the alphabet $A_2$) where
$\alpha=\la_e$, $\alpha+\beta=R_2$, and $\beta^+ = \la_w$.
In particular, $\la_w$ and $\la_e$ are complementary in
$R_2$.  It is not hard to see that this implies that
the shapes $\la_{ne}$ and $\la_{sw}$ are complementary in
$R_1\cap R_2$.

The above argument shows that if $\LRT(\la;(R_1,R_2))$
is nonempty then it is a singleton, since the entire tableau
$T\in \LRT(\la;(R_1,R_2))$ was specified.

For the converse, suppose $\la$ satisfies the two properties
given above.  Let $\la_w$ (resp. $\la_w$) be the part of
$\la-R_1$ in (resp. not in) the first $\mu_1$ columns.
Let $\alpha=\la_e$, $\beta=R_2-\alpha$, and 
$T$ the (not necessarily column-strict) tableau of shape $\la$
whose restrictions to the subshapes $R_1$, $\la_e$,
and $\la_w$ of $\la$ are given by $Y_1$, $\key(\alpha)$,
and $\key(\beta)$ where the key tableaux are taken with
respect to the alphabet $A_2$.  The column-reading word of $T$
satisfies $\cword(T)|_{A_1}\Knuth Y_1$ and
$\cword(T)|_{A_2}\Knuth Y_2$, so it only remains to show that
$T$ is indeed column-strict.  Since all the letters of $A_1$
are strictly smaller than those of $A_2$, the only possible violations
of column-strictness in $T$ are of the form
$T(r,\mu_1)>T(r,\mu_1+1)$ where the cell $(r,\mu_1)$ is in $\la_w$
and the cell $(r,\mu_1+1)$ is in $\la_e$.  Now $(r,\mu_1)\in\la_w$
(and $\la_w\subset R_2$ as partitions) implies that $r>\eta_1$ and
$\mu_2\ge \mu_1$, while $(r,\mu_1+1\in\la_e$ (and $\la_e\subset R_2$
as partitions) implies that $\eta_2\ge r$.  Now
\begin{equation*}
  T(r,\mu_1) = \key(\beta)(r-\eta_1,\mu_1)
\end{equation*}
is the $(r-\eta_1)$-th smallest letter in the last column of $\key(\beta)$,
which consists of the largest $\eta_2-m$ letters in
$A_2$, where $m$ is the number of parts of $\alpha$ equal to
$\mu_2$.  It follows that $T(r,\mu_1)=m+r$.  On the other hand,
\begin{equation*}
  T(r,\mu_1+1)=\key(\alpha)(r,1)=\eta_1+r,
\end{equation*}
being the $r$-th smallest letter in $A_2$.  By assumption
$m+r=T(r,\mu_1)> T(r,\mu_1+1)=\eta_1+r$, that is, $m>\eta_1$.
Now the partition $\alpha$ contains the cell $(m,\mu_2)$,
and $\la_e=\alpha$ as partitions, so it follows that
$\la$ contains the cell $(m,\mu_1+\mu_2)$.  But this cell lies outside
$R_1\cup R_2\cup\NE(R)\cup\SW(R)$ since $m>\eta_1=\min(\eta_1,\eta_2)$
and $\mu_1+\mu_2>\max(\mu_1,\mu_2)$.
\end{proof}

\begin{ex} Let $\mu_1=3$, $\mu_2=5$, $\eta_1=2$, $\eta_2=3$, and 
$\la=(76521)$.  Then in the notation of Proposition \ref{two rect},
$R_1\cup R_2=R_2=(555)$, $\la_{ne}=(210)$, $\la_{sw}=(210)$,
$\alpha=\la_e=(432)$ and $\beta=(123)$ so that $\la_w=(321)$.
The unique tableaux in $\LRT(\la;(R_1,R_2))$ and
$\LRT(\la;(R_2,R_1))$ are given respectively by
\begin{equation*}
  \begin{matrix}
     1&1&1&3&3&3&3\\
     2&2&2&4&4&4& \\
     3&4&5&5&5& & \\
     4&5& & & & & \\
     5& & & & & &
  \end{matrix}\qquad
  \begin{matrix}
     1&1&1&1&1&4&4\\
     2&2&2&2&2&5& \\
     3&3&3&3&3& & \\
     4&5& & & & & \\
     5& & & & & &
  \end{matrix}
\end{equation*}
\end{ex}

\subsection{Proof of Theorem \ref{action}}

\begin{proof} All will be proven here except the part of (A5)
given by the braid relation $\tau_p \tau_{p+1} \tau_p=
\tau_{p+1}\tau_p\tau_{p+1}$, which is proven later using the
induction coming from $R$-cocyclage.

For the well-definedness of the map $\tau_p$,
the only step requiring proof is (4), and this follows from (A1).

(A1): Since $w$ and $\tau_p w$ agree at all positions except those
occupied by letters of $B=A_p\cup A_{p+1}$, without loss of generality
it may be assumed (by restricting to $B$) that $R=(R_1,R_2)$ and $p=1$.
In this case, $Q(w)=Q(\tau_p w)$ by definition.  But a theorem
of D. White \cite{W} asserts that a word fits a skew shape $D$
if and only if its recording tableau satisfies a condition that depends
only on $D$.  Thus $w$ fits $D$ if and only if $\tau_p w$ does, proving (A1).

(A2): Without loss of generality suppose
$v\Knuth w$ is an elementary Knuth equivalence of words in $W(R)$.
Suppose first that $v=tyxzu$ and $w=tyzxu$ where
$t$ and $u$ are words of length $l$ and $m$ respectively
and $x\le y < z$ are letters.  This case is settled
by the following lemma.

\begin{lem} There exist words $t'$ and $u'$ of lengths $l$ and $m$
respectively, and letters $x'\le y'<z'$ such that
\begin{equation*}
  \tau_p v = t' y' x' z' u' \qquad \text{ and }\qquad
  \tau_p w = t' y' z' x' u'
\end{equation*}
\end{lem}
\begin{proof} Let $B=A_p\cup A_{p+1}=A_p'\cup A_{p+1}'$.
Suppose $x\not\in B$.  Then the removal of $x$ makes $v$ and $w$ identical;
call this common subword $v'$.  
By the definition of $\tau_p$, $\tau_p v$ and $\tau_p w$ are identical if the
letter $x$ is removed; this common subword is $\tau_p v'$.
It follows that $\tau_p v$ and $\tau_p w$ have the desired form.

The same argument works if $z\not\in B$.  So it may assumed that
$x\in B$ and $z\in B$.  Since $B$ is an interval and $x\le y<z$,
$y\in B$ as well.  Let $I$ be the set of positions of letters of $v$
that are not in $B$.  By the definition of $\tau_p$, the positions of
the letters not in $B$ are equal to $I$ for all of the words $v$,
$w$, $\tau_p v$ and $\tau_p w$, and all of these words agree at those
positions.  Thus by restricting to $B$, it may be assumed
that $R=(R_p,R_{p+1})$.  In this case we
have $Q(\tau_p v)=Q(v)$ and $Q(\tau_p w)=Q(w)$, so two applications of Lemma
\ref{Knuth} show that the words $\tau_p v$ and $\tau_p w$ have the
desired form.
\end{proof}

The other kind of elementary Knuth equivalence takes the form $v\Knuth w$
where $v=txzyu$ and $w=tzxyu$ with $t$ and $u$ as above and $x<y\le z$.  
Let $R\#=(R_t,R_{t-1},\dotsc,R_1)$ and $A_1\#$ through $A_t\#$ the
corresponding subalphabets for $R\#$.  Let $v\#$ be 
the reverse of the complement of the word $v$ in the alphabet $[n]$.
Then $v\# = u' y' z' x' t'$ and $w\# = u' y' x' z' t'$ where $t'$
and $u'$ are words of the same length as $t$ and $u$ respectively, and
$z'\le y'<x'$.  By Lemma \ref{rect ev}, $v\#$ and $w\#$ are in $W(R\#)$.
This puts us in the previous case, so by an application of the Lemma
we have $\tau_{t-p} v\# \Knuth \tau_{t-p} w\#$.  By Lemma \ref{rect ev},
$(\tau_p v)\# \Knuth (\tau_p w)\#$, which implies that
$\tau_p v\Knuth \tau_p w$ as desired.  This proves (A2).

(A3) follows immediately from (A2).

It remains to show (A4).  Suppose first that $p<t-1$.  Let
$A'=A_1\cup A_2\cup\dotsb\cup A_{t-1}$.  By definition the positions
of the letters of $A_t$ are the same in $w$ and $\tau_p w$.
By $\eta_t$ applications of Lemma \ref{hat} it is enough to show that
\begin{equation*}
  Q(w|_{A'}) = Q((\tau_p w)|_{A'})
\end{equation*}
since one uses the same process to pass from these tableaux to
$Q(w)$ and $Q(\tau_p w)$ respectively.  But this holds by induction
on $t$ since the latter tableau is equal to $Q(\tau_p(w|_{A'}))$.

Thus it may be assumed that $p=t-1$.  For the case $t=2$, (A4) holds by
definition.  So it may be assumed that $t>2$.
Let $R\#=(R_t,R_{t-1},\dotsc,R_1)$.  Then $Q(w)^{ev}=Q(w\#)$ and 
\begin{equation*}
  Q(\tau_p w)^{ev}=Q((\tau_p w)\#)=Q(\tau_1(w\#))
\end{equation*}
by Lemma \ref{rect ev}.  Since $\tau_1$ switches the first two rectangles
and $t>2$, $Q(\tau_1(w\#))=Q(w\#)$ by a previous argument.
Thus $Q(w)^{ev}=Q(\tau_p w)^{ev}$ and $Q(w)=Q(\tau_p w)$,
proving (A4).

(A5):  The fact that $\tau_p$ is an involution, and that
$\tau_p \tau_q = \tau_q \tau_p$ if $|p-q|>1$, follow easily from the
definitions.

(A6): This follows from the special case of the operator $\tau_p$ 
and interval $I=[p,p+1]$.

(A7): This immediately reduces to the case (1) where $R=(R_1,R_1)$.
But $\tau_1:\LRT(\la;(R_1,R_1))\rightarrow \LRT(\la;(R_1,R_1))$
is the identity since $\LRT(\la;(R_1,R_1))$ is empty or a singleton.
\end{proof}

\subsection{Proof of Proposition \ref{chi}}

\begin{rem} \label{chi obs} Let $w$ be a word of content
$\gamma(R)=(\mu_1^{\eta_1},\dotsc,\mu_t^{\eta_t})$,
the content of any word in $W(R)$.
\begin{enumerate}
\item For each $1\le j\le t$, the set of positions occupied by
the letters of $A_j$ in $\chi_R^k(w)$ are rotated cyclically to the right
$k$ positions from the corresponding set of positions in $w$ (by
Theorem \ref{action} (A6) applied to the automorphism of
conjugation $w_0^R$).
\item If $j\not=i$, then the letters of $A_j$ just move to the right by one
position (and are otherwise unchanged) in passing from $w$ to $\chi_R(w)$
(by 1 and Theorem \ref{action}(A7) applied to the automorphism of conjugation
$w_0^R$ and the word $w$, which has $\mu_j$ copies of each letter in
the interval $A_j$).
\item Let $B=A_i\cup A_{i+1}\cup \dotsb A_j$ and $w=uv$ where $v$ has
length $k$.  Say $v|_B$ has length $k'$.  Then
\begin{equation*}
  \chi_R^k(uv)|_B = \chi_R^{k'}(u|_B\, v|_B).
\end{equation*}
This follows from the case $k=1$, where 2 applies.
\end{enumerate}
\end{rem}

\begin{proof} Let $w\in W(R)$.  For 1, it suffices to prove that
$\chi_R(w)\in W(R)$, since 2 implies that $\chi_R^N$ is the identity,
where $N$ is the length of the word $w$ (making $\chi_R$ invertible).
Write $w=ux$ with $x\in A_i$ say.
For $j\not=i$, $\chi_R(w)|_{A_j}=w|_{A_j} \Knuth Y_j$ by
Remark \ref{chi obs} and the fact that $w\in W(R)$.
For $j=i$, $\chi_R(w)|_{A_i} = \chi_R(w|_{A_i})$.  Thus
it may be assumed that $R=(R_1)$.  Let $R_1=(a^n)$, $Y=Y_1$, and
$w_0^R=w_0$.  For this it certainly suffices to let $uv\in W((R_1))$ and
show that $(w_0 v)(w_0 u)\in W(R)$, for we only need the case that
$v$ is a letter.  Let $\nu$ be the content of $v$.
By Lemma \ref{one rect}, $\nu$ is a partition, $P(v)=\key(\nu)$,
and $u\Knuth \key((a^n)-\nu)$.  Then by Lemma \ref{one rect}
and the definition of key tableau,
\begin{equation*}
\begin{split}
  P((w_0 v)(w_0 u)) &= P((w_0 \key(\nu))(w_0 \key((a^n)-\nu))) \\
  &=P(\key(\rev(\nu))\key((a^n)-\rev(\nu))) \\
  &=\key((a^n))
\end{split}
\end{equation*}
where $\rev(a_1,\dotsc,a_n)=(a_n,\dotsc,a_1)$.

Note that 2 holds by definition when $k=1$.
In light of Remark \ref{chi obs} 1, it is enough to show that
\begin{equation*}
  \chi_R^k(uv)|_{A_i} = ((w_0^R v)(w_0^R u))|_{A_i}
\end{equation*}
for all $i$.  Fix $i$.  Let $k'$ be the length of $v|_{A_i}$.  Then
by Remark \ref{chi obs} 3,
\begin{equation*}
  \chi_R^k(uv)|_{A_i} = \chi_R^{k'}((u|_{A_i} v|_{A_i})).
\end{equation*}
Thus we have reduced to the case that $R=(R_1)$.  Let $R_1=(a^n)$,
$Y=Y_1$, and $w_0=w_0^R$.  By the bijectivity of the RS correspondence
it is enough to show that $\chi_R^k(uv)$ and $(w_0 v)(w_0 u)$ have
the same $P$ tableaux and the same $Q$ tableaux.  By 1 and its proof,
both words have $P$ tableau equal to $Y_1$.

Let $Q=Q(w)=Q(uv)$, $Q'=Q((w_0 v)(w_0 u))$ and
$Q"=Q(\chi_R^k(uv))$.  It only remains to show that $Q'=Q"$.
Recall that all of these tableaux have shape equal to that
of $Y$, which is rectangular with $n$ columns and $a$ rows.
Let $N=k+l$ be the number of cells in $Q$.
In light of Lemma \ref{rect ex}, it is enough to show
that $Q' = \pr_k(Q)$ and $Q" = \pr_1^k(Q)$.  Let $T+j$ denote
the tableau whose entries are obtained from those of $T$ by
adding the integer $j$.  We have
\begin{equation*}
  P(Q|_{[l+1,l+k]})-l = Q(v) = Q(w_0 v) = Q'|_{[k]}
\end{equation*}
by Lemma \ref{jeu recording}, Theorem \ref{action} for $w_0$,
and the definition of recording tableau.  Also
\begin{equation*}
  P(Q'|_{[k+1,k+l]})-k=Q(w_0 u)=Q(u)= Q|_{[l]}.
\end{equation*}
It follows that $Q' = \pr_k(Q)$.  To show that $Q"=\pr_1^k(Q)$,
it suffices to show that $Q(\chi_R(w))=\pr_1(Q)$.  But the
above argument with $k=1$ proves this.
\end{proof}

\subsection{Proof of Theorem \ref{cocyc action}}

The proof requires a few preliminary results.

\begin{lem} \label{chi Knuth} Suppose $u y x z\in W(R)$ with $x<y\le z$
letters, so that $u y x z \Knuth u y z x$ is an elementary
Knuth equivalence.  Then $\chi_R^3(u y x z) \Knuth \chi_R^3(u y z x)$.
\end{lem}
\begin{proof} We have
\begin{equation*}
  \chi_R^3(u y x z) = (w_0^R yxz)(w_0^R u) \Knuth
  (w_0^R yzx)(w_0^R u) = \chi_R^3(u y z x)
\end{equation*}
by two applications of Proposition \ref{chi} and Theorem \ref{action}.
\end{proof}

Before giving the proof of Theorem \ref{cocyc action} it is useful
to state a more detailed version in the two rectangle case.
Let $R=(R_1,R_2)$ and $R'=s_1 R=(R_2,R_1)$.  Suppose
$w\in W(R)$, written $w=ux$ with $x$ a letter.  Let $\la$
(resp. $\rho$) be the shape of $P(w)=P(ux)$ (resp.
$P(\chi_R(w))=P((w_0^R x)(w_0^R u))$).
Let $s$ (resp. $s'$) be the cell giving the difference of $\la$
(resp. $\rho$) and  the shape of $P(u)$, which is the same as 
the shape as $P(w_0^R u)$.

\begin{prop} \label{two cyc action} With the above notation
and that of Proposition \ref{two rect}:
\begin{enumerate}
\item If $s\in\NE(R)$ then $s'\in \SW(R)$, $\rho_{ne}=\la_{ne}-\{s\}$,
and $x\in A_2$.
\item If $s\in\SW(R)$ then $s'\in \NE(R)$, $\rho_{se}=\la_{se}-\{s\}$,
and $x\in A_1$.
\item If $s\in R_1\cup R_2$ then $s'=s$ and $\rho=\la$.
If $s\in R_1$ and $s\not\in R_2$ then $x\in A_1$.
If $s\in R_2$ and $s\not\in R_1$ then $x\in A_2$.
(The corner cell $s$ cannot lie in $R_1\cap R_2$).
\end{enumerate}
In particular, $\chi_R$ and $\tau_1$ commute for $R=(R_1,R_2)$.
\end{prop}

The proof, which relies on explicit computations
on two-rectangle LR tableaux, is straightforward but tedious
and is omitted.  In each of the three cases there are two subcases
depending on whether one rectangle contains the other or not.

Proof of Theorem \ref{cocyc action}:
\begin{proof} Let $A_1'$, $A_2'$, etc. be the subalphabets for $R'$
and $Y_1'$, $Y_2'$, etc. the Yamanouchi tableaux.
Let $w=ux\in W(R)$ with $x$ a letter.  Say $x\in A_i$.

Suppose $i\not\in\{p,p+1\}$.  We have
\begin{equation*}
  \chi_{R'}(\tau_p w) = \chi_{R'} (\tau_p u) x =
  (w_0^{R'} x)(w_0^{R'} \tau_p u)
\end{equation*}
\begin{equation*}
  \tau_p \chi_R w = \tau_p (w_0^R x)(w_0^R u) = (w_0^R x)(\tau_p w_0^R u)
\end{equation*}
since $x\in A_i$ and $w_0^R x\in A_i$ are unchanged by $\tau_p$.
But $w_0^{R'}x=w_0^R x$ so it enough to show that
$w_0^{R'}\tau_p u=\tau_p w_0^R u$.  It is clear that these two
words agree at positions of letters not in $B=A_p\cup A_{p+1}$,
so it is enough to show that
$w_0^{R'}\tau_p (u|_B)=\tau_p w_0^R (u|_B)$.
But $w_0^R u|_B = u|_B$ and $w_0^{R'} \tau_p (u|_B)=\tau_p u|_B$
by Theorem \ref{action}(A7) for the automorphisms of conjugation
$w_0^R$ and $w_0^{R'}$.  Thus both words are equal to $\tau_p u|_B$.

The other case is that $i\in\{p,p+1\}$.  For each $j\not\in\{p,p+1\}$
it follows from the definitions and Remark \ref{chi obs} that
$\chi_{R'}(\tau_p w)$ and $\tau_p (\chi_R(w))$ agree at the positions
containing letters of $A_j=A_j'$.  Again it is enough to show that
\begin{equation*}
  \chi_{R'}(w\tau_p)|_B = ((\chi_R(w))\tau_p)|_B
\end{equation*}
where $B=A_p \cup A_{p+1} = A_p'\cup A_{p+1}'$.  Let $\tau_p w = u' x'$
where $x'$ is a letter.  Since $w=ux$ with $x\in B$, $x'\in B$
by the definition of $\tau_p$.  We have
\begin{equation*}
\begin{split}
  \chi_{R'}(\tau_p w)|_B &= \chi_{R'}((\tau_p w)|_B) =
  \chi_{R'}(\tau_p(w|_B)) \\
  (\tau_p(\chi_R(w)))|_B &= \tau_p((\chi_R(w))|_B)
	\tau_p(\chi_R(w|_B))
\end{split}
\end{equation*}
by Remark \ref{chi obs} and the definition of $\tau_p$.

Thus we have reduced to the two rectangle case $(R_p,R_{p+1})$.
For simplicity of notation let $t=2$, $p=1$, and $\tau=\tau_1$.
By the bijectivity of the RS correspondence it
suffices to show that the two words
$\chi_{R'}(\tau w)$ and $\tau \chi_R(w)$ have the same $P$ and $Q$ tableaux.
For this we claim that it suffices to show that their $P$ tableaux have the
same shape.  Indeed, since we are in the two rectangle case and
both $P$ tableaux are in the set $\LRT(R)$ and are assumed to have the
same shape, they must be equal by Proposition \ref{two rect}.
For the equality of the $Q$ tableaux, it suffices to show that they agree
after applying the invertible operator $\pr_1$ (see section \ref{rec tab}).
But since it is assumed that these $Q$ tableaux have the same shape,
it is enough to show that
\begin{equation} \label{defl Q}
  P(Q(\chi_{R'}(\tau w))|_{[2,N]}) =
  P(Q(\tau(\chi_R(w)))|_{[2,N]}).
\end{equation}
Let $w=ux$ and $\tau w=u'x'$ where $x$ and $x'$ are letters.
Then $\chi_R(w)=(w_0^R x)(w_0^R u)$
and $\chi_{R'}(\tau w)=(w_0^R x')(w_0^R u')$.
Let $N$ be the length of $w$.  We have
\begin{equation*}
\begin{split}
  P(Q(\tau(\chi_R(w)))|_{[2,N]})-1
  &= P(Q(\chi_R(w))|_{[2,N]}) - 1 \\
  &= P(Q((w_0^R x)(w_0^R u))|_{[2,N]}) -1 \\
  &= Q(w_0^R u) = Q(u) = Q(w)|_{[N-1]}
\end{split}
\end{equation*}
by Theorem \ref{action}, definition of $\chi_R(w)$, Lemma
\ref{jeu recording}, Theorem \ref{action}, and the
fact that $w=ux$.  On the other hand,
\begin{equation*}
\begin{split}
  P(Q(\chi_{R'}(\tau w))|_{[2,N]}) - 1
  &= P(Q((w_0^{R'} x')(w_0^{R'} u'))|_{[2,N]}) -1 \\
  &= Q(w_0^{R'}u') = Q(u') = Q(\tau w)|_{[N-1]} = Q(w)|_{[N-1]}
\end{split}
\end{equation*}
for similar reasons.  But this establishes \eqref{defl Q}.

So it only remains to show that the shapes of
$W:=P(\chi_R(w))$ and $W':=P(\chi_{R'}(\tau w))$ coincide,
since $W$ has the same shape as $\tau W=P(\tau(\chi_R(w)))$
by Theorem \ref{action}.  Let $s$ be the cell giving the difference
of the shapes of $P(w)$ and $P(u)$; this is also the difference
of the shapes of $P(w\tau)$ and $P(u')$, by Theorem \ref{action}
and the previously proven fact that $Q(u)=Q(u')$.
Propositions \ref{two cyc action} and \ref{two rect}
explicitly show how the shape of $P(w)$ (resp. $P(\tau w)$),
together with the cell $s$, determines the shape of
$W$ (resp. $W'$).  But $P(w)$ and $P(\tau w)$ have the same shape,
so $W$ and $W'$ do as well.
\end{proof}

\subsection{Proof of Theorem \ref{poset thm}}

\begin{proof} First it is shown that $\LRT(R)$ under $\le_R$
is a partial order.  It is enough to show that
$\le_R$ has an extension by a partial order $\le$.
Let $n_j(T)$ be the number of letters in $T|_{A_j}$ in the first $\mu_j$
columns.  Define the partial order $\le$ on $\LRT(R)$ by $T < S$ if there
is an index $j$ such that $n_j(T)=n_j(S)$ for all $j<i$ but
$n_i(T) > n_i(S)$.  Suppose $T <_R S$ is a covering relation
coming from the word $w=ux\in W(R)$ where $x$ is a letter.  Let
$x\in A_i$ say.  $x$ is the smallest letter in $A_i$
since $x$ is the last letter of $w|_{A_i}$ and $P(w|_{A_i})=Y_i$, by 
Lemma \ref{one rect}.

Let $s$ be the cell giving the difference of the shapes of $U:=P(u)$
and $S=P(w)$.  Write $\chi_R(ux)=x'u'$ so that $T=P(x'u')$.
Write $U':=P(u')=w_0^R U$.  Now $x'=w_0^R x$ is the largest
letter of $A_i$ and $u'=w_0^R u$.  For $j<i$ we have
\begin{equation*}
  n_j(T) = n_j(U) = n_j(U') = n_j(S).
\end{equation*}
The first equality holds, for by a property of Schensted insertion,
since the cell $s$ is not in the first $a$ columns,
the tableaux $S$ and $U$ must agree in the first $a$ (and hence $\mu_i$)
columns.  The second equality comes from Theorem \ref{action} (A7)
for the automorphism of conjugation $w_0^R$.  The third
follows from the fact that $x'\in A_i$ is strictly greater
than any letter in $A_j$, so the column insertion of $x'$ into
$U'$ doesn't move any letters in $A_j$.  Moreover,
\begin{equation*}
  n_i(T) = n_i(U) = n_i(U') = n_i(S)-1.
\end{equation*}
The first two equalities hold for the above reasons.
Recall that $S$ is obtained from $U'$ by the column insertion
of the letter $x'$.  But $x'$, being the maximum letter of
$A_i$, only bumps other copies of $x'$ in $U'$.
There are $\mu_i-1$ copies of $x'$ in $U'$, hence at most that
many in the first $\mu_i$ columns of $U'$.  Let $B=A_1\cup\dotsb\cup A_i$.
Then by explicit calculation $P((x' U')|_B)$ is obtained from $U'|_B$ by
adjoining the letter $x'$ to the bottom of the first column
of $U'|_B$ that does not contain the letter $x'$.
This puts another letter of $A_i$ into the first $\mu_i$ columns.
Thus $S<T$.

This proves that $\LRT(R)$ is a poset under $\le_R$.

Before showing that $\LRT(R)$ is graded, let us prove 2 and 3.
Let $T\in \LRT(R)$.  If $T$ has more than $a$ columns, then it
has a corner cell in a column strictly east of the $a$-th,
and therefore admits a covering relation $T' <_R T$.  Continuing this
process one produces a saturated chain in $\le_R$ from $T$ down to an
element $\min(T)\in\LRT(R)$ that has at most $a$ columns.
Since $\min(T)|_{A_i}\Knuth Y_i$ and $a=\mu_i$ for some $i$,
$\min(T)$ must have exactly $a$ columns.  For 3, note
that for any $T\in\LRT(R)$, $T|_{A_1}=Y_1$.  For $T$ minimal
$T$ must consist of $Y_1$ atop a tableau $\Th$.  By Remark
\ref{rect rem}, $\Th\in\LRT(\Rhat)$.

To show $\LRT(R)$ is graded, it is enough to show that
for any $T\in\LRT(R)$, there is only one such tableau $\min(T)$
(that is, a tableau in $\LRT(R)$ that has $a$ columns and satisfies
$\min(T)\le_R T$), and any saturated chain from $T$
down to $\min(T)$ has the same length.  Let $T\in\LRT(R)$
and $S_1$ and $S_2$ in $\LRT(R)$ with $a$ columns and suppose that
there are saturated chains from $T$ down to $S_1$ and $S_2$.
Let $T_1<_R T$ and $T_2<_R T$ be the first steps
in these saturated chains from $T$ down to $S_1$ and $S_2$ respectively.
We construct a tableau $T_3$ that is two steps below both $T_1$
and $T_2$.  To see that this suffices, by induction and definition
one has $\min(T_3)=\min(T_1)=S_1$ and $\min(T_3)=\min(T_2)=S_2$,
so $S_1=S_2$ and $\min(T)$ is the common $\min$ of
$T$, $T_1$, $T_2$, and $T_3$.  Furthermore the
distance from $T_1$ and $T_2$ down to the common minimum
is well-defined and equal, since both distances are two more than
the distance from $T_3$ down to the common minimum.
So it suffices to construct $T_3$.
Its construction uses a rectangular analogue of \cite[Lemme 2.13]{LS2}.

Let $s_1$ and $s_2$ be the two corner cells of the shape of
$T\in\LRT(\la;R)$, both strictly east of the $a$-th column,
which induce the covering relations $T_1 <_R T$ and $T_2 <_R T$.
Without loss of generality assume that
$s_1$ is strictly north and strictly east of $s_2$.
Let $(U_i,x_i)$ be the pair consisting of a column-strict tableau
$U_i$ of shape  of shape $\la-\{s_i\}$ and a letter $x_i$
for $1\le i\le 2$, computed by reverse row insertion on $T$
at $s_i$.  Clearly there is a corner cell $s_3$ of the shape
$\la-\{s_1,s_2\}$ that is strictly east of the $a$-th column, such
that $s_2$ is strictly south and weakly west of $s_3$ and
$s_1$ is strictly east and weakly north of $s_3$.  Performing
reverse row insertions on $T$ at the cells $s_1$, then $s_2$, then
$s_3$, let $V_1$ be the resulting column-strict tableau
and $yxz$ the three ejected letters.  Doing the same thing
except using the order $s_2$, then $s_1$, then $s_3$, let
$V_2$ be the resulting tableau and $y'z'x'$ the three ejected letters.
By Lemma \ref{Knuth} we have $V_1=V_2$ and $y'=y$, $z'=z$ and $x'=x$
with $x\le y<z$, and $T \Knuth U' y x z \Knuth U' y z x$.
Let $T_3=P(\chi_R^3(U'yxz))=P(\chi_R^3(U'yzx))$; the latter equality
holds by Lemma \ref{chi Knuth}.  There are saturated chains
\begin{equation*}
\begin{split}
  &P((\chi_R^3(U'yxz))<_R P((\chi_R^2(U'yxz))<_R
  	P((\chi_R(U'yxz))=T_1<_R P(U'yxz) \\
  &P((\chi_R^3(U'yzx))<_R P((\chi_R^2(U'yzx))<_R
  	P((\chi_R(U'yzx))=T_2<_R P(U'yzx),
\end{split}
\end{equation*}
where both left hand tableaux are $T_3$ and both right hand tableaux are
$T$.  The fact that these covering relations are produced by corner cells
that are strictly east of the $a$-th column, is a consequence of
Lemma \ref{SE cell}.
\end{proof}

\subsection{Proof of Proposition \ref{switch iso}}
\begin{proof} $\tau_p$ is a bijection so it is enough
to show that if $T <_R S$ is a covering relation,
then $\tau_p T <_{R'} \tau_p S$ is a covering relation.

Let $T\in \LRT(\la;R)$, $s$ a corner cell of $\la$
that is strictly east of the $a$-th column where $a:=\max_i \mu_i$,
$U$ the column-strict tableau of shape $\la-\{s\}$ and $x$ the letter
such that $T=P(Ux)$, and $W=P(\chi_R(Ux))$.  Then $W <_R T$,
and all covering relations have this form.

Let $\tau=\tau_p$ and $N$ be the number of cells in $T$.
Let $U'x'=\tau(Ux)$ and $T'=\tau T$.
By Theorem \ref{action}, $T$ and $T'$ have the same shape
and $Q(U')=Q(U'x')|_{[N-1]}=Q(Ux)|_{[N-1]}=Q(U)$, so that
$U'$ is a column-strict tableau of the same shape as $U$.
Define $W':=P(\chi_{R'}(U'x'))$.  By Theorem \ref{cocyc action},
$W'=\tau W$.  But by construction, $W' <_R T'$.
\end{proof}

With this poset structure in place, we now show that
$\tau_p \tau_{p+1} \tau_p =\tau_{p+1}\tau_p \tau_{p+1}$,
the only part of Theorem \ref{action}(A5) that was not proven.
\begin{proof} Let $B=A_p\cup A_{p+1}\cup A_{p+2}$.  Consider the
operators $\tau_p \tau_{p+1} \tau_p$ and $\tau_{p+1}\tau_p \tau_{p+1}$.
Both do not disturb the letters outside the interval $B$, so by restriction
to $B$, we may assume that $R=(R_1,R_2,R_3)$ and $p=1$.  By
Theorem \ref{action}(A3) and (A4) and the bijectivity of the
RS correspondence, it is enough to check the equality of the two
operators on $T\in\LRT(\la;R)$.  Using the fact that $\tau_1$ and
$\tau_2$ are involutions, we may reorder $R$ so that $\mu_1$
is the largest among the $\mu_i$.  Proceeding by induction on height
in the poset $\LRT(R)$ with order $\le_R$, by Theorem
\ref{poset thm} it may be assumed that $T$ is minimal,
that is, $T$ has exactly $\mu_1$ columns.  Using Remark \ref{narrow action}
below, it is easy to see that $\tau_1$ and $\tau_2$ satisfy the
braid relation on $T$.
\end{proof}

\begin{rem} \label{narrow action} Let $T\in\LRT(\la;R)$ where
$\la$ has $a=\max_i \mu_i$ columns with $\mu_j=a$.  Let
$A_i'$ be the $i$-th subalphabet for $\tau_j R$.  Then $\tau_j T$
is obtained from $T$ by a rather trivial vertical exchanging process.
In each column of $T$, first replace the letters of $A_j$ (resp. $A_{j+1}$)
by their counterparts in $A_{j+1}'$ (resp. $A_j'$)
and then sort the resulting column; only the letters
in $A_j\cup A_{j+1}=A'_j\cup A'_{j+1}$ need to be moved.
The key fact is that each column of $T|_{A_j}$ must contain each of the
numbers in $A_j$.
\end{rem}

\subsection{Proof of Proposition \ref{chi charge}}

\begin{proof} Observe that
\begin{equation} \label{d and tau}
d_{i,R}(w)=d_{i,\tau_i R}(\tau_i w)
\end{equation}
since $\tau_i$ preserves shape by Theorem \ref{action}.
Consider first the case that $t=2$, that is, $R=(R_1,R_2)$.
Let $R'= \tau_1 R$.  Then
\begin{equation*}
\begin{split}
  2! (\LLT_R(w) - \LLT_R(\chi_R(w))) &=
  	d_R(w)-d_R(\chi_R(w)) +
	d_{R'}(\tau_1 w)-d_{R'}(\chi_R(\tau_1 w)) \\
	&= 2(d_R(w)-d_R(\chi_R(w)))
\end{split}
\end{equation*}
by definition, the fact that $\chi$ and $\tau_1$ commute
(Theorem \ref{cocyc action}), and \eqref{d and tau}.
The values of $d_R(w)$ and $d_R(\chi_R(w))$ are determined by
the shapes of the two-rectangle LR tableaux $P(w)$ and $P(\chi_R(w))$.
By Lemma \ref{two cyc action}, it follows that
\begin{equation} \label{two diff}
  d_R(w)-d_R(\chi_R(w))=
  \begin{cases}
    1 & \text{if $f(id)=f(\tau_1)=2$} \\
    -1 & \text{if $f(id)=f(\tau_1)=1$} \\
    0 & \text{otherwise.}
  \end{cases}
\end{equation}

In general we have
\begin{equation} \label{mult LLT diff}
  t! (\LLT_R(w) - \LLT_R(\chi_R(w))) =
  \sum_{\tau\in S_t} \sum_{i=1}^{t-1} (t-i)
  (d_{i,\tau R}(\tau w) - d_{i,\tau w}(\chi_R(\tau w))).
\end{equation}
Next it is shown that
\begin{equation} \label{d diff}
  d_{i,\tau R}(\tau w) - d_{i,\tau w}(\chi_R(\tau w)) =
  \begin{cases}
    1 & \text{if $f(\tau)=f(\tau \tau_i)=i+1$} \\
    -1 & \text{if $f(\tau)=f(\tau \tau_i)=i$} \\
    0 & \text{otherwise.}
  \end{cases}
\end{equation}
Suppose first that $f(\tau)\not\in\{i,i+1\}$, that is,
the last letter of $\tau w$ is not in the union $B$ of the
$i$-th or $(i+1)$-st subalphabets for $\tau R$.  Then
$(\chi_R(\tau w))|_B=(\tau w)|_B$ and the difference on the
left hand side of \eqref{d diff} is zero.  If
$f(\tau)\in\{i,i+1\}$, then one can deduce
\eqref{d diff} using the formula \eqref{two diff} for the
two-rectangle case.  This proves \eqref{d diff}.

As by-product of the above calculation we obtain
\begin{equation} \label{f val}
  f(\tau_i \tau) \in \{ f(\tau), \tau_i f(\tau) \}
\end{equation}
for any $\tau\in S_t$ and $1\le i\le t-1$.

For a statement $P$ let $\chi(P)$ be $1$ if $P$ holds and $0$ otherwise.
Define
\begin{equation*}
  \chargediff(f) := \sum_{\tau\in S_t} \sum_{i=1}^{t-1}
  (t-i) (\chi(f(\tau)=f(\tau\tau_i)=i+1)-
         \chi(f(\tau)=f(\tau\tau_i)=i)).
\end{equation*}
By \eqref{d diff} we have
\begin{equation*}
\chargediff(f) = t!(\LLT_R(w)-\LLT_R(\chi_R(w))).
\end{equation*}

To compute $\chargediff(f)$, it is convenient to replace $f$
by a simpler function $f'$ such that $\chargediff(f')=\chargediff(f)$.
Suppose $f$ has maximum value $M+1$ for $1<M<t$.  Define
\begin{equation*}
  f'(\tau)=\begin{cases}
    f(\tau) & \text{if $f(\tau)\not=M+1$} \\
    M & \text{if $f(\tau)=M+1$.}
   \end{cases}
\end{equation*}
By definition one has the equality of inverse images
\begin{equation} \label{inverse image}
  f^{-1}(i) = (f')^{-1}(i) \qquad\text{if $i\not\in\{M,M+1\}$.}
\end{equation}
We have
\begin{equation} \label{cdiff}
\begin{split}
&\quad\chargediff(f)-\chargediff(f') \\
&=\sum_{i=1}^{t-1} (t-i) \sum_{\tau\in S_t}
    (\chi(f(\tau)=f(\tau_i\tau)=i+1)
	-\chi(f'(\tau)=f'(\tau_i\tau)=i+1) \\
	&-\chi(f(\tau)=f(\tau_i\tau)=i)
	+\chi(f'(\tau)=f'(\tau_i\tau)=i) 
    ).
\end{split}
\end{equation}
The summand is zero unless $\{i,i+1\}\cap\{M,M+1\}\not=\emptyset$ or
equivalently $i\in\{M-1,M,M+1\}$.  For $i\in[t-1]$, let $Y_{i,\tau}$
be the inner sum of \eqref{cdiff}.  We have
\begin{equation*}
\begin{split}
  Y_{M-1,\tau} &= \chi(f(\tau)=f(\tau_{M-1}\tau)=M)-
  	\chi(f'(\tau)=f'(\tau_{M-1}\tau)=M) \\
  	&= \chi(f(\tau)=f(\tau_{M-1}\tau)=M)-
  		\chi(\{f(\tau)=f(\tau_{M-1}\tau\in\{M,M+1\})) \\
  	&= -\chi(f(\tau)=M+1)
\end{split}
\end{equation*}
by \eqref{inverse image}, the definition of $f'$ and \eqref{f val},
and \eqref{f val} again.  We have
\begin{equation*}
\begin{split}
Y_{M,\tau} &= \chi(f(\tau)=f(\tau_M\tau)=M+1)
-\chi(f(\tau)=f(\tau_M\tau)=M) \\
&+\chi(f'(\tau)=f'(\tau_M\tau)=M) \\
&= \chi(f(\tau)=f(\tau_M\tau)=M+1)
-\chi(f(\tau)=f(\tau_M\tau)=M) \\
&+\chi(\{f(\tau),f(\tau_M\tau)\}\subseteq\{M,M+1\}) \\
&= 2 \chi(f(\tau)=f(\tau_M\tau)=M+1) +
\chi(\{f(\tau),f(\tau_M\tau)\}=\{M,M+1\})
\end{split}
\end{equation*}
by the fact that $f'$ has maximum value $M$,
the definition of $f'$, and considering the four cases for values
of $f(\tau)$ and $f(\tau_M\tau)$ in the set $\{M,M+1\}$.  We have
\begin{equation*}
\begin{split}
Y_{M+1,\tau} &= -\chi(f(\tau)=f(\tau_{M+1}\tau)=M+1) \\
&= -\chi(f(\tau)=M+1))
\end{split}
\end{equation*}
since $f$ has maximum value $M+1$ and $f'$ has maximum value $M$,
and \eqref{f val}.
Define
\begin{equation*}
\begin{split}
  S_0 &= \{\tau\in S_t | f(\tau)=f(\tau_M\tau)=M+1\} \\
  S_1 &= \{\tau\in S_t | f(\tau)=M+1 \text{ and }f(\tau_M\tau)=M\}.
\end{split}
\end{equation*}
Note that $f^{-1}(M+1)$ is the disjoint union of
$S_0$ and $S_1$, and there is a bijection between $S_1$ and the set
\begin{equation*}
 \{\tau\in S_t | f(\tau)=M \text{ and }f(\tau_M\tau)=M+1\}.
\end{equation*}
Putting the above calculations together we have
\begin{equation*}
\begin{split}
\chargediff(f) - \chargediff(f') &=
  -(t-(M-1))|f^{-1}(M+1)|
  +(t-M)(2 |S_0| + 2 |S_1|) \\
  &-(t-(M+1)) |f^{-1}(M+1)| = 0.
\end{split}
\end{equation*}
In this calculation, the only property of $f$ that
was used was \eqref{f val}.  Furthermore it is easy to check that
$f'$ satisfies \eqref{f val} and $(f')^{-1}(1)=f^{-1}(1)$.
Repeating this construction, we obtain a function $h:S_t\rightarrow[t]$
that has maximum value 2 and satisfies $\chargediff(h)=\chargediff(f)$,
\eqref{f val}, and $h^{-1}(1)=f^{-1}(1)$.
Let us explicitly calculate $\chargediff(h)$.
Let $n_{ij}$ be the number of $\tau\in S_t$ such that
$h(\tau)=i$ and $h(\tau_1\tau)=j$.  Then we have
\begin{equation*}
\begin{split}
  t! &= n_{11} + n_{12} + n_{21} + n_{22} \\
  n_{21} &= n_{12} \\
  n_{11}+n_{12} &= |h^{-1}(1)| = |f^{-1}(1)| \\
  n_{21}+n_{22} &= |h^{-1}(2)| = t!-|f^{-1}(1)|
\end{split}
\end{equation*}
Since $h$ has maximum value 2,
the condition $h(\tau)=h(\tau_2\tau)=2$ is equivalent to
$h(\tau)=2$.  Then 
\begin{equation*}
\begin{split}
  \chargediff(f)
  &= \chargediff(h) \\
  &= -(t-1)n_{11} + (t-1)n_{22} - (t-2)(n_{21}+n_{22}) \\
  &= -(t-1)n_{11} -(t-2)n_{12} + n_{22} \\
  &= -(t-1)n_{11} -(t-2)n_{12} + (t! - n_{11} - 2 n_{12}) \\
  &= t! - t\, n_{11} - t\, n_{12} \\
  &= t! - t\, |f^{-1}(1)| \\
  &= t! - t!\, |\NC|.
\end{split}
\end{equation*}
The last equality follows from the fact that for $1\le i\le t$
the elements $\tau_1 \tau_2\dotsb \tau_{i-1}$ form a system of
coset representatives of $S_1\times S_{t-1}\backslash S_t$,
$f^{-1}(1)$ is a union of the corresponding cosets, each of which
has cardinality $(t-1)!$.
\end{proof}

\subsection{Proof of Theorem \ref{explicit charge}}
	
\begin{proof} To show that $\LLT_R$ is constant on Knuth classes
it is enough to show that $w\mapsto d_{i,\sigma R}(\sigma w)$ is,
where $\sigma$ acts on $w$ by a composite of $\tau_p$'s.  Write
$R'_1$ and $R'_2$ for $(\sigma R)_i$ and $(\sigma R)_{i+1}$
and $B$ the union of the $i$-th and $(i+1)$-st alphabets 
for $\sigma R$.  Recall that
\begin{equation*}
  d_{i,\sigma R}(\sigma w)=d_{R'_1,R'_2}((\sigma w)|_{A'_1\cup A'_2}).
\end{equation*}
Now $v\Knuth w$ implies $\sigma v\Knuth \sigma w$ which implies
$(\sigma v)|_B\Knuth (\sigma w)|_B$ by Theorem \ref{action} (A2)
and the fact that $B$ is an interval.  But then
$P((\sigma v)|_B)=P((\sigma w)|_B)$, and the map
$d_{i,\sigma R}$ depends only on this $P$ tableau.

(C1) certainly holds for $\LLT_R$.  For (C2), in light
of Proposition \ref{chi charge}, it is enough to show that the set
$\NC$ is empty.  Write $T' <_R T$ where $w=ux\in W(R)$,
$T=P(w)$, $T'=P(\chi_R(w))$, and $s$ is the cell where the shapes
of $T=P(w)$ and $U=P(u)$ differ.  Since $s$ lies in a column strictly
east of the $a$-th (where $a=\max_i \mu_i$), it follows that
$x$ came from a cell of $T$ in a column strictly east of the $a$-th.
This means $x\not\in A_1$, since $T|_{A_1}=Y_1$ lies entirely in the
first $\mu_1$ columns, and $\mu_1\le a$.  The same argument shows that
for any permutation $\sigma \in S_t$, the reverse row insertion on
$\sigma T$ at $s$ ejects a number that cannot be in the first subalphabet
of $\sigma R$.  This shows $\NC$ is empty, so that (C2) holds.

(C3): For $\tau\in S_t$, define $\sigma\in S_{t-1}$ by
\begin{equation*}
  \sigma(i) = \begin{cases}
    \tau(i+1) & \text{if $\tau(i+1)<\tau(1)$} \\
    \tau(i+1)-1 & \text{if $\tau(i+1)>\tau(1)$.} \\
  \end{cases}
\end{equation*}
There is a bijection $\tau\leftrightarrow (\sigma,\tau(1))$.
By Remark \ref{narrow action} and Theorem \ref{two rect},
\begin{equation*}
  d_{i,\tau R}(\tau T) = \begin{cases}
    d_{i,\sigma\Rhat}(\sigma\Th) & \text{if $1\le i<\tau(1)$} \\
    0 & \text{if $i\in\{\tau(1),\tau(1)+1\}$}\\
    d_{i-1,\sigma\Rhat}(\sigma\Th) & \text{if $\tau(1)+1<i\le t-1$}
  \end{cases}
\end{equation*}
Then letting $j=\tau(1)$ we have
\begin{equation*}
\begin{split}
  t!\, \LLT_R(T) &=  \sum_{\tau\in S_t}
	\sum_{i=1}^{t-1} (t-i) \,d_{i,\tau R}(\tau T) \\
  &= \sum_{j=1}^t \sum_{\sigma\in S_{t-1}}
    (\sum_{i=1}^{j-1} (t-i)\, d_{i,\sigma\Rhat}(\sigma\Th)\\
  &+ \sum_{i=j+2}^{t-1} (t-i)\, d_{i-1,\sigma\Rhat}(\sigma\Th)) \\
  &= \sum_{\sigma\in S_{t-1}} \sum_{j=1}^t
	( \sum_{i=1}^{j-1} (t-i)\, d_{i,\sigma\Rhat}(\sigma\Th) \\
  &+ \sum_{i=j+1}^{t-2} (t-i-1)\, d_{i,\sigma\Rhat}(\sigma\Th) ) \\
  &= \sum_{\sigma\in S_{t-1}} \sum_{j=1}^t
     (\sum_{i=1}^{j-1}\, d_{i,\sigma\Rhat}(\sigma\Th) \\
	&-(t-j-1)\,d_{j,\sigma\Rhat}(\sigma\Th)+
	\sum_{i=1}^{t-2}(t-i-1)\, d_{i,\sigma\Rhat}(\sigma\Th)) \\
  &= t!\, \LLT_{\Rhat}(\Th) +
     \sum_{\sigma\in S_{t-1}}
	(-\sum_{j=1}^{t-1} (t-j-1)\,d_{j,\sigma\Rhat}(\sigma\Th)) \\
	&+\sum_{1\le i<j\le t-1}\, d_{i,\sigma\Rhat}(\sigma\Th)) \\
  &= t!\, \LLT_{\Rhat}(\Th) +
     \sum_{\sigma\in S_{t-1}}
	(-\sum_{j=1}^{t-1} (t-j-1)\,d_{j,\sigma\Rhat}(\sigma\Th)) \\
	&+\sum_{1\le i\le t-2} (t-i-1)\, d_{i,\sigma\Rhat}(\sigma\Th)) \\
  &=  t!\, \LLT_{\Rhat}(\Th).
\end{split}
\end{equation*}

(C4) holds by the definition of $\LLT_R$.
\end{proof}

\subsection{Proof of Theorem \ref{intrinsic charge_R}}

\begin{proof} For existence, one has the explicit function
$\LLT_R$ by Theorem \ref{explicit charge}.

For uniqueness, the proof proceeds by induction on $\LLT_R$,
then on the number of inversions of $R$ (the number of pairs
$1\le i<j\le t$ such that either $\mu_i<\mu_j$ or
$\mu_i=\mu_j$ and $\eta_i<\eta_j$), then on the number $t$ of
rectangles in $R$.

Suppose $T$ is not $<_R$-minimal.  Then by Remark \ref{corner chi}, (C2)
applies and drops $\LLT_R$.  Otherwise suppose $T$ is $<_R$-minimal.

Now suppose $R$ has an inversion.  Then it has an adjacent inversion,
say $(p,p+1)$.  Then (C4) applies and $\tau_p R$ has fewer inversions
than $R$ does.  Otherwise suppose $R$ has no inversions.

In this case $\mu_1\ge \mu_2\ge \dotsb\ge\mu_t$, so
(C3) applies and drops $t$ by one.
\end{proof}

\subsection{Proofs of Lemmas for main theorem}

Proof of Lemma \ref{preserving map}:
\begin{proof} Let $(w,T,U)\in\SSS$ and $\Phi(w,T,U)=(w,P,Q)$.
First it must be shown that $(w,P,Q)\in\TTT$, that is,
the word of $P$ is $R$-LR.  
By direct computation we have 
\begin{equation*}
  y^n \dotsb y^1 Y_1 \Knuth
  y^n \dotsb y^{m+1} m^{\mu_1} y^m (m-1)^{\mu_1} y^{m-1} \dotsb 1^{\mu_1} y^1.
\end{equation*}
Let $d=|\alpha(w)|-|R_1|$.  Since none of the cycled letters are in the
first alphabet $A_1=[m]$, we have
\begin{equation*}
\begin{split}
  \word(T)\word(Y_1) &\Knuth u^m \dotsb u^1 u^n \dotsb u^{m+1} \word(Y_1) \\
  &= \chi_R^d(y^n \dotsb y^{m+1} \word(Y_1) y^m \dotsb y^1) \\
  &\Knuth\chi_R^d(y^n \dotsb y^{m+1} m^{\mu_1} y^m (m-1)^{\mu_1} y^{m-1}
	\dotsb 1^{\mu_1} y^1 )\\
  &=\chi_R^d(v^n \dotsb v^1) \Knuth \chi_R^d(P),
\end{split}
\end{equation*}
which holds by explicit Knuth equivalences and Lemma \ref{chi Knuth}.

By assumption $\word(T)$ is $\Rhat$-LR, so by definition
$\word(T)\word(Y_1)$ is $R$-LR.
But both Knuth equivalence and $\chi_R$ preserve $R$-LRness,
so $\word(P)$ is $R$-LR and $\Phi$ is well-defined.

By definition $\Phi$ is injective, being a composition of injective maps.
$\Phi$ is sign-preserving by definition.

To show that $\Phi$ is weight-preserving, it suffices to show that
each of the $d$ instances of the operator $\chi_R$ in the above
computation, induce a $\le_R$-covering relation.  That is, if 
\begin{equation*}
  W_i := P(\chi_R^i(y^n \dotsb y^{m+1} \word(Y_1) y^m \dotsb y^1))
\end{equation*}
then it must be shown that $W_{i+1}<_R W_i$.  Let us say that
a cell is sufficiently east if it lies in a column strictly east of
the $\mu_1$-st.

Fix $0\le i\le d-1$.  Write
\begin{equation*}
  u x v = y^n \dotsb y^{m+1} \word(Y_1) y^m \dotsb y^1
\end{equation*}
where $v$ has length $i$ and $x$ is a letter.
Also write $w_0^R ux = u' x'$ and $w_0^R v = v'$.
By Theorem \ref{chi} $W_i=\chi_R^i(uxv)=v' u' x'$.  It is enough to
show that the cell $s$ given by the difference of the shapes of
$P(v'u'x')$ and $P(v'u')$ is sufficiently east.  
By Proposition \ref{SE cell} it is enough to
show that the cell $s'$ given by the difference of the shapes of $P(u'x')$
and $P(u')$ is sufficiently east.  The cell $s'$ is also the
difference of the shapes of $Q(u'x')$ and $Q(u')$.  Now
$Q(u'x')=Q(w_0^R ux)=Q(ux)$ by Theorem \ref{action}(A4).
By the definition of recording tableau it follows that
$s'$ is the difference of the shapes of $Q(ux)$ and $Q(u)$, or
equivalently, of $P(ux)$ and $P(u)$.

Now $ux$ is an initial subword of $\word(Y_1) y^m \dotsb y^1$
and $ux$ contains $\word(Y_1)$.  Each of the words $y_i$ is a weakly
increasing word consisting of letters that are strictly greater
than those in the tableau $Y_1$.  By Theorem \ref{Pieri} it follows
that the tableau $P(Y_1 y^m \dotsb y^1)$ has at most $m$ rows
and consists of the tableaux $Y_1$ and $P(y^m \dotsb y^1)$ sitting
side by side.  In particular the cells of the difference
of the shapes $P(Y_1 y^m\dotsb y^1)$ and $Y_1$ is sufficiently east,
which implies that $s'$ is sufficiently east.
\end{proof}

Proof of Lemma \ref{inv lem}:
\begin{proof} By definition $\theta$ is sign-reversing and weight-
preserving on all of $\TTT$.  By Lemma \ref{lattice and pairs} $\theta$ is
an involution.  It remains to show that $\theta$ stabilizes the set
$\Phi(\SSS)$.  Let $(w,P,Q)\in\Phi(\SSS)$ and $\theta(w,P,Q)=(w',P',Q')$.
It may be assumed that $(w,P,Q)$ is not a fixed point of $\theta$.
Let $v'$ be to $(P',Q')$ as $v$ is to $(P,Q)$ in the definition of
$\Phi$.  It is enough to show that $(v')^{i}$ starts with the subword
$i^{\mu_1}$ for every $1\le i\le m$, since the other steps in the map
$\Phi$ are invertible by definition.

Let us apply Lemma \ref{two row dual} (see section \ref{two row stuff})
to $Q$ and $Q' = s_r e_r Q$.
There is nothing to prove unless $r\le m$.
Suppose first that $r<m$.  We need only check that
${v'}^r$ starts with $r^{\mu_1}$ and
${v'}^{r+1}$ starts with $(r+1)^{\mu_1}$.
Since $v^r = r^{\mu_1} y^r$ and
$v^{r+1} = (r+1)^{\mu_1} y^{r+1}$ where
all the letters of $y^r$ and $y^{r+1}$ are strictly greater than $m$,
it follows that
\begin{equation*}
  P({v'}^{r+1} {v'}^r)|_{[r,r+1]} =
  P(v^{r+1} v^r)|_{[r,r+1]} = P((r+1)^{\mu_1} r^{\mu_1}).
\end{equation*}
This, together with the fact that ${v'}^{r+1}$ and ${v'}^r$ are weakly
increasing words, implies that all of the letters $r+1$ must precede all
of the letters $r$ in the word ${v'}^{r+1} {v'}^r$, that is,
${v'}^i$ starts with $i^{\mu_1}$ for $i\in[r,r+1]$.

The remaining case is $r=m$.  Then $v^r=r^{\mu_1} y^r$
and $v^{r+1} = y^{r+1}$.  Let us calculate ${v'}^{r+1}$ and ${v'}^r$
using a two-row jeu-de-taquin.  Let $V$ (resp. $V'$) be the (skew) two row
tableau with first row $v^r$ (resp. ${v'}^r$) and second row $v^{r+1}$
(resp. ${v'}^{r+1}$) in which the two rows achieve the maximum overlap.
The overlaps of $V$ and $V'$ are equal by Lemma \ref{overlap} and
the fact that $Q$ and $Q'$ are in the same $r$-string
and hence have the same $r$-paired letters.  Furthermore this
common overlap is at least $\mu_1$.  To see this, note that
the overlap weakly exceeds the minimum of $\mu_1$ and $|y^{r+1}|$
since all of the letters in $y^{r+1}$ have values in the alphabet
$[m+1,n]=[r+1,n]$ and there are $\mu_1$ copies of $r$ in $v^r$.
On the other hand, $|y^{r+1}|>\mu_1$, for otherwise by
Lemma \ref{overlap} all of the letters $r+1$ in $Q$ would be
$r$-paired, contradicting the choice of $r$.

We calculate $V'$ from $V$ in two stages.
Let $V"$ be the two row skew tableau (whose rows have maximum
overlap) such that $P(V")=P(V)$, where the first row of $V"$ is
one cell longer than that of $V$.  By Lemma \ref{two row dual} this
tableau exists since $Q$ has an $r$-unpaired letter $r+1$; the
corresponding recording tableau is $e_r Q$.  Furthermore $V"$ is
obtained by sliding the ``hole" in the cell just to the left of the first
letter in the first row of $V$, into the second row.
By the same reasoning as above, $V"$ has the same overlap that $V$ does.
Finally we calculate $V'$ from $V"$ by another two row jeu-de-taquin.
If the first row of $V"$ is shorter than the second, we are done,
for in this case the first row ${v'}^r$ of $V'$ contains the first
row of $V"$, which in turn contains the first row of $V$, which
contains $r^{\mu_1}$.  So suppose
the second row of $V"$ is shorter than the first, by $p$ cells, say.
Now $p$ is less than or equal to the number of cells on the right
end of the first column of $V"$ that have no cell of $V"$ below them.
Since $V"$ has maximum overlap it follows that
when $p$ holes are slid from the second row of $V"$ to the first,
they all exchange with numbers lying in the portion of the first row
of $V"$ that extends properly to the right of the second.
Thus the subword $r^{\gamma_r}$ remains in the first row of $V'$,
and we are done.
\end{proof}

\begin{ex} In $v'$ all the subwords are the same as in $v$ except
that ${v'}^2 = 2225677$ and ${v'}^3=333567$.  In this example $r=m$.
The tableaux $V$, $V"$, and $V'$ are given below.
\begin{equation*}
\begin{split}
V&=\begin{matrix}
  \sq&\sq&\sq&2&2&2&5&6 \\
  3&3&3&5&6&7&7&7
\end{matrix} \\
R"&=\begin{matrix}
  \sq&\sq&2&2&2&5&6&7 \\
  3&3&3&5&6&7&7&
\end{matrix} \\
R'&=\begin{matrix}
  \sq&2&2&2&5&6&7&7 \\
  3&3&3&5&6&7& &
\end{matrix}
\end{split}
\end{equation*}
\end{ex}

\section{Schensted miscellany}
\label{Schensted}

This section is the repository for some well-known (or should be
well-known) facts regarding the Robinson-Schensted-Knuth correspondence.

\subsection{Evacuation}

Given a word $w=w_1\dotsb w_N$ in the alphabet $[n]$, let
$w\#_n$ be the reverse of the complement of $w$ with respect to the
alphabet $[n]$, that is, the $i$-th letter of $w\#_n$ is
$n+1-w_{N+1-i}$ for $1\le i\le N$.

Given a column-strict tableau $T$ of partition shape in the
alphabet $[n]$, define $T^{ev_n}$ to be the 
unique column-strict tableau in the alphabet $[n]$ such that
the shape of $(T^{ev_n})|_{[i]}$ is equal to
that of $P(T|_{[n+1-i,n]})$ for all $1\le i\le n$.

\begin{thm} \label{ev} \cite{LS2} Let $w$ be a word of length $N$ in the
alphabet $[n]$, $P=P(w)$ and $Q=Q(w)$.  Then $P(w\#_n)=P^{ev_n}$ and
$Q(w\#_n)=Q^{ev_N}$.
\end{thm}

Applying evacuation to LR tableaux produces other LR tableaux.

\begin{lem} \label{rect ev} Let $R=(R_1,\dotsc,R_t)$ and
$R\#=(R_t,R_{t-1},\dotsc,R_1)$.  Then $\#_n$ and $ev_n$ restrict to
bijections such that the following diagram commutes:
\begin{equation*}
\begin{CD}
  W(R) @>{RS}>> \bigcup_\la \LRT(\la;R) \times \ST(\la) \\
  @V{\#_n}VV 			@VV{\cup_\la (ev_n\times ev_N)}V \\
  W(R\#) @>>{RS}> \bigcup_\la \LRT(\la;R\#) \times \ST(\la)  
\end{CD}
\end{equation*}
\end{lem}
\begin{proof} Let $A_1'$ through $A_t'$ be the subalphabets
for $R\#$ and $Y_1'$ through $Y_t'$ the corresponding Yamanouchi tableaux.
Let $w\in W(R)$, $\#=\#_n$ and $ev=ev_n$.  In light of Theorem \ref{ev}
it is enough to check that $w\#$ is in $W(R\#)$.
\begin{equation*}
  P((w\#)|_{A_p'}) = P((w|_{A_p})\#) =
  P(w|_{A_p})^{ev} = (Y_p)^{ev} = Y_p'.
\end{equation*}
Thus $w\#\in W(R\#)$.
\end{proof}

\subsection{Removal of large letters}

\begin{lem} \label{hat} Suppose $w=w_1 w_2\dotsb w_N$ has maximum letter
$n$, occurring in positions $1\le i_1<i_2<\dotsb<i_r\le N$,
and let $\wh$ be the word obtained by removing these letters $n$ from
$w$.  Let $P=P(w)$, $Q=Q(w)$, $\Ph=P(\wh)$ and $\Qh=Q(\wh)$ where
$\wh$ is recorded in $\Qh$ by the letters in $[N]-\{i_1,\dotsc,i_r\}$.
Then $\Ph = P(w)|_{[n-1]}$ and $Q$ is obtained from $\Qh$ by the
row insertion of the letters $i_1$ through $i_r$ in that order.
\end{lem}

\subsection{Pieri's rules}
\begin{thm} \label{Pieri} Let $w=w_1 w_2\dotsb w_N$ and $Q=Q(w)$
its row-insertion recording tableau.  Then
\begin{enumerate}
\item If $w_i\le w_{i+1}$ then $i+1$ is strictly east and weakly
north of $i$ in $Q$.
\item If $w_i>w_{i+1}$ then $i+1$ is strictly south and weakly west of
$i$ in $Q$.
\end{enumerate}
\end{thm}

\subsection{Knuth equivalence and recording tableaux}

\begin{lem} \label{Knuth} Let $v$ and $w$ be words with
$P(v)=P(w)$.  The following are equivalent.
\begin{enumerate}
\item $v\Knuth w$ is an elementary Knuth equivalence of the form
\begin{equation*}
  v = t y x z u \qquad\text{ and }\qquad
  w = t y z x u
\end{equation*}
where $t$ and $u$ are words $ x\le y<z$ are letters with $x$ and $z$
in positions $i$ and $i+1$ in $v$.
\item $Q(v)$ and $Q(w)$ differ by the transposition of $i$ and $i+1$,
and in $Q(v)$, $i$ is strictly south and weakly west of $i-1$ and
$i+1$ is strictly east and weakly north of $i-1$.
\end{enumerate}
\end{lem}

\subsection{Recording tableaux}
\label{rec tab}

The following result is not hard to prove
and appears as \cite[Lemma 2.4.4]{Sh}.

\begin{lem} \label{jeu recording}
Let $w=uv$ be a word where $u$ and $v$ have lengths
$N-k$ and $k$ respectively.  Let $Q=Q(uv)$.  Then
\begin{equation*}
  P(Q|_{[N-k+1,N]}) = Q(v)+N-k
\end{equation*}
where $Q(v)+N-k$ is the tableau formed by adding the number $N-k$
to every entry in $Q(v)$.
\end{lem}

Let $Q$ be a standard tableau with $N$ letters.
For the integer $i$, define $\pr_i(Q)$ to be the standard tableau
such that
\begin{equation*}
\begin{split}
  \pr_i(Q)|_{[N-i]} &= P(Q|_{[i+1,N]})-i \\
  P(\pr_i(Q)|_{[N-i+1,N]} &= Q|_{[i]}+N-i
\end{split}
\end{equation*}
$\pr_i(Q)$ can be defined in terms of jeu-de-taquin or
exchanging tableaux \cite{BSS}; it is obtained by exchanging the
subtableaux $Q|_{[1,i]}$ and $Q_{[i+1,N]}$ (and then relabeling).
For standard tableaux $\pr_1$ is Sch\"utzenberger's promotion
operator \cite{Schu}.

The following lemma is not hard to prove using the
jeu-de-taquin techniques of \cite{H}.  It is crucial that the shape of
$Q$ be both normal (having a unique northwestmost cell)
and antinormal (having a unique southwestmost cell), that is,
the shape of $Q$ must be rectangular.

\begin{lem} \label{rect ex}
Let $Q$ be a standard tableau of rectangular shape.
Then $\pr_k(Q) =\pr_1^k(Q)$.
\end{lem}

\begin{lem} \label{SE cell}
Let $U$ be a column-strict tableau of partition shape,
$x$ a letter, and $v$ a word.  Let $s$ (resp. $s'$)
be the cell given by the difference of the shapes of $P(Ux)$
(resp. $P(vUx)$) and $U$ (resp. $P(vU)$).  Then $s'$ is
weakly south and weakly east of $s$.
\end{lem}
\begin{proof} Let $Q$ be the skew standard tableau that records
the column insertion of $v$ into $P(Ux)$.  Then $s'$ is
the vacated cell after a jeu-de-taquin that slides $Q$ to the
northwest into the cell $s$.  This precise statement follows
from \cite[Lemma 21]{Sh}.  More crudely, the tableau
$P(vU)$ is entrywise smaller than $U$, viewing empty cells
as containing the letter $\infty$.  So the row insertion of
$x$ into $P(vU)$ must necessarily end at a cell $s'$
that is weakly south and weakly west of the cell $s$
where the row insertion of $x$ into $U$ ends.
\end{proof}

\subsection{Two row jeux-de-taquin}
\label{two row stuff}

There is a duality between the crystal operators and
jeux-de-taquin on two-row skew column strict tableaux.
This is described below.

Define the \textit{overlap} of the pair $(v,u)$ of weakly increasing
words to be the length of the second row in the tableau $P(v u)$,
or equivalently, the maximum number of columns of size two among the
skew column strict tableaux with first row $u$ and second row $v$.

The following results appear in \cite{SW}.

\begin{lem} \label{overlap} Let $(P,Q)$ be the tableau pair
obtained by column RSK from the sequence of words
$\{v^i\}$.  Then the overlap of the pair of words
$(v^{r+1},v^r)$ is equal to the number of $r$-pairs in $Q$.
\end{lem}

\begin{lem} \label{two row dual} Let $\{v^i\}$ and $(P,Q)$ be
as in Lemma \ref{overlap} and let $\{{v'}^i\}$ be another
sequence of weakly increasing words with corresponding tableau
pair $(P',Q')$.  The following are equivalent.
\begin{enumerate}
\item $P=P'$, and $Q$ and $Q'$ are in the same $r$-string.
\item ${v'}^i = v^i$ for $i\not\in[r,r+1]$ and
$P({v'}^{r+1} {v'}^r) = P(v^{r+1} v^r)$.
\end{enumerate}
\end{lem}

\end{document}